\documentclass[12pt]{article}
\usepackage{graphicx}
\usepackage{amsmath}
\usepackage{amssymb}
\usepackage{graphicx}
\usepackage{graphics}
\usepackage{amsthm}
\usepackage{amscd}
\usepackage{color}
\usepackage{amsfonts}
\usepackage{fancyhdr}
\usepackage{graphicx}
\usepackage{latexsym,amssymb}
\usepackage{amsthm}
\usepackage{indentfirst}
\usepackage{amsmath}
\usepackage{color}
\usepackage{xcolor}
\usepackage{fourier}
\usepackage[all]{xy}
\usepackage[colorlinks=true]{hyperref}
\usepackage{authblk}

\newtheorem{theoremalph}{Theorem}

\newtheorem*{Main Theorem}{Main Theorem}
\newtheorem{Coro}[theoremalph]{Corollary}
\newtheorem{Theorem}{Theorem}[section]
\newtheorem*{Theorem A}{Theorem A}
\newtheorem*{Theorem A'}{Theorem A'}
\newtheorem*{Theorem B'}{Theorem B'}

\newtheorem{Definition}[Theorem]{Definition}
\newtheorem{Proposition}[Theorem]{Proposition}
\newtheorem{Lemma}[Theorem]{Lemma}

\newtheorem{Remark}[Theorem]{Remark}

\newtheorem{Corollary}[Theorem]{Corollary}

\newtheorem{Claim-numbered}[Theorem]{Claim}

 \def\NN{{\mathbb N}} 

 \def\RR{{\mathbb R}}

 \def\ZZ{{\mathbb Z}}

\def\La{\Lambda}

\def\cA{{\cal A}}   \def\cM{{\cal M}} 
    
   \def\cO{{\cal O}} \def\cU{{\cal U}}
    \def\cV{{\cal V}}
    \def\cW{{\cal W}}
\def\cF{{\cal F}}

\newcommand{\diff}{{\operatorname{Diff}}}

\def\diff{\operatorname{Diff}}

\def\dim{\operatorname{dim}}

\def\ud{\operatorname{d}}
\def\e{{\varepsilon}}

\def\det{\operatorname{det}}

\def\Leb{\operatorname{Leb}}
\def\homeo{\operatorname{Homeo}}

\def\m{\operatorname{m}}
\begin{document}

\title{Entropy properties of   mostly expanding partially hyperbolic  diffeomorphisms}
\author{Jinhua Zhang}


\maketitle

\begin{abstract}
	The statistical properties of mostly expanding partially hyperbolic diffeomorphisms have been substantially studied. 
	In this paper, we would like to address the entropy properties of mostly expanding partially hyperbolic diffeomorphisms.  
	We prove that 
  for  mostly expanding partially hyperbolic diffeomorphisms with minimal strong stable foliation and one-dimensional center bundle, there exists  a $C^1$-open neighborhood of them,  in which the topological entropy varies continuously and the intermediate entropy property holds. To prove that, we  show that each non-hyperbolic ergodic measure is approached by horseshoes in entropy and in weak$*$-topology.
 
\end{abstract}
%
%
%

\tableofcontents

\section{Introduction}
 Entropy (topological entropy and metric entropy) is an important topological invariant and plays a central role in describing the complexity of a dynamical system. The dependence of entropy on the systems is an interesting topic and has been substantially studied by many mathematicians. 
  Smoothness of a system is  relevant to the upper semi-continuity of entropy. For $C^\infty$-diffeomorphisms, Newhouse \cite{N} proved that the metric entropy varies upper semi-continuously and  Yomdin \cite{Yo} proved that the topological entropy also varies upper semi-continuously. For systems with lower regularity, a bifurcation phenomena called \emph{homoclinic tangency} (i.e. the existence of non-transverse homoclinic intersections between the stable and unstable manifolds of a hyperbolic saddle) seems to be an obstruction to the upper semi-continuity of topological entropy. Misiurewicz \cite{Mi} constructed $C^r$-diffeomorphisms with homoclinic tangencies at which the entropy fails to be upper semi-continuous. For $C^1$-diffeomorphisms far away from homoclinic tangencies, one can recover the upper semi-continuity of entropy (see \cite{LVY,DFPV}). 
   As for the lower semi-continuity of the topological entropy, the hyperbolicity is involved in.  
 For $C^1$-uniformly  hyperbolic diffeomorphisms, the topological entropy is locally constant due to Smale's structural stability theorem. Using Pesin's theory, Katok \cite{Ka} proved that hyperbolic ergodic measures of  $C^{1+\alpha}$ surface diffeomorphisms are approached by hyperbolic horseshoes and thus the topological entropy is lower semi-continuous  for $C^{1+\alpha}$ surface diffeomorphisms. As a consequence, the topological entropy varies continuously among $C^\infty$ surface diffeomorphisms.

 Then it is natural to ask:  among which classes of differentiable systems, does the topological entropy vary continuously?
 The examples in \cite{Mi} and the results in \cite{LVY,DFPV} tell us that one  should consider the systems away from homoclinic tangencies. Such diffeomorphisms   are partially hyperbolic with multi-one dimensional center bundles due to \cite{CSY}. Therefore,   it is reasonable to firstly consider the continuity of topological entropy for partially hyperbolic diffeomorphisms with one-dimensional center bundle.  
There are three types of classical examples of partially hyperbolic diffeomorphisms with one dimensional center bundle:
  skew-products with circle fiber,
  derived from Anosov diffeomorphisms and 
   time one maps of Anosov flows. 
The topological entropy of these classical examples  varies continuously. See for instance \cite{BFSV,U,SY}. In recent years, some new anomalous examples appeared. The continuity of the topological entropy of the anomalous  examples in \cite{BGP, BGHP} was proved in \cite{YZ}.  Saghin and Yang \cite{SY} conjectured  that \emph{the topological entropy varies continuously among the set of partially hyperbolic diffeomorphisms with one-dimensional center bundle}. 

 The classical variational principle gives the relation between metric entropy and topological entropy: topological entropy equals the supremum of the metric entropy.  For partially hyperbolic diffeomorphisms with one-dimensional center bundle,   the usual approach to prove the continuity of the topological entropy is to find horseshoes approximating ergodic measures in weak$*$-topology and in entropy. This makes the problem of the  continuity of topological  entropy closely related to the intermediate entropy property, since   horseshoes satisfy the intermediate entropy property.  


 Recall that a diffeomorphism $f$ satisfies  the \emph{intermediate entropy property} if for any $h\in[0,h_{top}(f))$, there exists an ergodic measure $\nu$ such that $h_\nu(f)=h.$ In general, one can not expect that $h$ is chosen as $h_{top}(f)$, since the ergodic measures whose metric entropies obtain the topological entropy do not always exist \cite{Mi}.  
  Katok proposed a conjecture that 
\emph{$C^r$ ($r\geq 1$)-diffeomorphisms satisfy the intermediate entropy property}.  It is classical that hyperbolic systems satisfy the intermediate entropy property. Katok's  result \cite{Ka} implies  that $C^r$ ($r>1$) non-uniformly hyperbolic systems have the intermediate entropy property. In general, for dynamics beyond uniform hyperbolicity, this conjecture is widely open and there are not so  many results. See for instance \cite{S3,S,S2,GSW,YZ,LSWW}. Even for partially hyperbolic diffeomorphisms with one-dimensional center bundle, Katok's conjecture is still open. See \cite{S3,S,YZ} for some partial results. 

There is a special class of partially hyperbolic diffeomorphisms called mostly expanding diffeomorphisms whose study was initiated in \cite{ABV}\footnote{The notion in \cite{ABV} nowadays is  called partially hyperbolic diffeomorphisms with non-uniformly expanding center and it is slightly different from the notion of mostly expanding that we use.  See the discussions and examples in \cite{AV}}. It was proved in  \cite{AV} that  mostly expanding partially hyperbolic diffeomorphisms form an open set.  There are many  works on studying such kind of diffeomorphisms,  and people aim to show the statistical properties of such systems, for instance the existence and finiteness of physical measures or SRB measures, statistical stability of the set of SRB measures and the mixing properties of SRB measures. See for instance \cite{A,AV,AV2,AlDLV,AlLi, Y} and references therein.  But there are few results on the entropy properties of such systems. In this paper, we are interested in studying the continuity of topological entropy and intermediate entropy properties of such systems.  The method we use is finding horseshoes approximating ergodic measures in weak$*$-topology and in entropy. 

\subsection{Statements of our main results}
A diffeomorphism $f\in\diff^1(M)$ is \emph{partially hyperbolic} if there exist a $Df$-invariant continuous splitting $TM=E^s\oplus E^c\oplus E^u$ and some constants $C>1,\lambda\in(0,1)$ such that for any $x\in M$ and $n\in\mathbb{N}$, one has
\begin{itemize}
	\item $\|Df^n|_{E^s(x)}\|<C\lambda^n \textrm{~~and~~}\|Df^{-n}|_{E^u(x)}\|<C\lambda^n.$
	\item $\|Df^n|_{E^s(x)}\|\cdot \|Df^{-n}|_{E^c(f^n(x))}\|<C\lambda^n \textrm{ and }\|Df^n|_{E^c(x)}\|\cdot \|Df^{-n}|_{E^u(f^n(x))}\|<C\lambda^n. $
	\end{itemize}
By \cite{Go}, up to changing a metric, one can assume that $C=1$, and we will do so throughout this paper.
Due to \cite{HPS}, the bundles $E^s$ and $E^u$ are uniquely integrable to $f$-invariant foliations, called strong stable and strong unstable foliations and denoted by $\cF^s$ and $\cF^u$ respectively. We denote by $\cF^*(x)$ the $\cF^*$-leaf through the point $x$, for $*=s,u.$  

Let $\mu$ be an invariant measure of   a partially hyperbolic diffeomorphism $f$ with one dimensional center bundle (i.e. $\dim(E^c)=1$), then we define 
\[\chi^c(\mu,f):=\int\log\|Df|_{E^c}\|\ud\mu.\]
When $\mu$ is ergodic, by Oseledec's theorem and Birkhoff's ergodic theorem, $\chi^c(\mu,f)$ coincides with the   Lyapunov exponent of $\mu$ along the center bundle $E^c$ (called \emph{center Lyapunov exponent}).
When there is no ambiguity, we will simply write $\chi^c(\mu)$. 
%

For a partially hyperbolic diffeomorphism $f$, we denote by $G^u(f)$ the set of invariant measures satisfying the entropy formula for unstable entropy (see Equation~\eqref{eq:unstable-entropy-formula}).  We remark that for $f\in\diff^{1+\alpha}(M)$, the set $G^u(f)$ coincides with the set of u-Gibbs states (see Theorem~\ref{thm.ledrappier}), and moreover, $f$ is called \emph{mostly expanding} if all the center Lyapunov exponents of each $u$-Gibbs state are positive (see Section~\ref{s.mostly-expanding} for more information).  

 Recall that a foliation is \emph{minimal} if every leaf is dense in the manifold. 
 \begin{theoremalph}\label{thm.existence}
 	Let $f\in\diff^{1}(M)$ be a  partially hyperbolic diffeomorphism with $\dim(E^c)=1$. 
 	Assume that 
 	\begin{itemize}
 			\item $\chi^c(\mu)>0$ for any $\mu\in G^u(f);$
 		\item  the strong stable foliation is minimal.
 	\end{itemize}  
 	Then there exists a $C^1$-open neighborhood $\mathcal{U}\subset \diff^1(M)$ of $f$ such that 
 	\begin{itemize}
 		\item[--]  each $g\in\cU$ satisfies the intermediate entropy property;
 		\item[--] the entropy function $g\in\cU\mapsto h_{top}(g)$ is continuous.
 	\end{itemize}
 	\end{theoremalph}
 \begin{Remark}
 	In fact,   each   $g\in\cU$ also satisfies our assumptions.
The first assumption is an open condition due to the upper semi-continuity of the compact set $G^u(f)$ (see Lemma~\ref{l.the-set-Gu(f)}).
 In general, the second assumption is not an open condition, but combining with the first assumption, one can deduce  that the strong stable foliation is $C^1$-robustly minimal (see Theorem~\ref{thm.robust-minimal-ss}). 
 \end{Remark}	
We can apply our result to some conservative partially hyperbolic diffeomorphisms, and obtain the continuity of topological entropy and  intermediate entropy property.  
 \begin{Coro}\label{cor:conservative-case}
 	Let $f\in\diff_{\m}^{1+\alpha}(M)$ be a    partially hyperbolic diffeomorphism preserving a smooth volume $\m$ and  with $\dim(E^c)=1.$ 
 	Assume that 
 	\begin{itemize}
 		\item   $\int\log\|Df|_{E^c}\|\ud \m>0$; 
 		\item  the strong stable foliation is minimal.
 	\end{itemize}  
 	Then there exists a $C^1$-neighborhood $\mathcal{U}\subset \diff^1(M)$ of $f$ such that 
 	\begin{itemize}
 		\item[--]  each $g\in\cU$ satisfies the intermediate entropy property;
 		\item[--] the entropy function $g\in\cU\mapsto h_{top}(g)$ is continuous.
 	\end{itemize}
 \end{Coro}
 \begin{Remark}
One can deduce that   $f$ is $C^1$-stably ergodic using    the minimality of strong stable foliation and the     center Lyapunov exponent  being positive (see Theorem~\ref{thm.citerion-mostly-expanding}).
 \end{Remark}
For the classical examples of partially hyperbolic diffeomorphisms (skew-products with circle fiber,
derived from Anosov diffeomorphisms and time one maps of Anosov flows), the results in \cite{BF,TY,CroPol} show that under some mild open conditions,  ergodic measures with high entropy are hyperbolic. However, for the systems that we  consider, it is not clear for us if there exist hyperbolic ergodic measures of maximal entropy (i.e. ergodic measures whose metric entropies equal the topological entropy),   and we are not able to exclude the existence of non-hyperbolic ergodic measures of maximal entropy.  To show the continuity of topological entropy and intermediate entropy property, we prove that non-hyperbolic ergodic measures are approached by hyperbolic horseshoes in entropy and in weak$*$-topology.

 Given an invariant compact set $K\subset M$ of $f\in\diff^1(M)$, we denote by $\mathcal{M}_{inv}(f,K)$ and $\mathcal{M}_{erg}(f,K)$ the sets of $f$-invariant   and $f$-ergodic measures on $K$ respectively. When $K=M$, we will simply write $\mathcal{M}_{inv}(f)$ and $\mathcal{M}_{erg}(f)$. Recall that an $f$-invariant compact set $\Lambda$ is called a \emph{hyperbolic basic set} if it is a hyperbolic transitive set and there exists a neighborhood $U$ of $\Lambda$ such that  $\cap_{n\in\ZZ}f^n(\overline{U})=\Lambda.$ 
  \begin{theoremalph}\label{thm.approaching-non-hyperbolic-ergodic-measure}
  		Let $f\in\diff^1(M)$ be a   partially hyperbolic diffeomorphism with $\dim(E^c)=1.$ 
  	Assume that 
  	\begin{itemize}
  				\item $\chi^c(\mu)>0$ for any $\mu\in G^u(f);$
  		\item  the strong stable foliation is minimal.
  	\end{itemize}  
  	Then there exist  a $C^1$-neighborhood  $\cU\subset\diff^1(M)$ of $f$,  and constants  $\kappa>0$ and $\chi_1>0$ such that for  any $g\in\cU$,   any ergodic measure $\nu\in\mathcal{M}_{\rm erg}(g)$ with $-\chi_1\leq \chi^c(\nu)\leq 0$ and any $\e>0$, there exists a hyperbolic basic set $\Lambda_\e$ of $g$ whose center bundle is uniformly expanding such that 
		\begin{itemize}
			\item[--] \[h_{top}(g, \Lambda_{\e})>\frac{h_\nu(g)-\e}{1+\kappa\cdot(|\chi^c(\nu)|+\e)}; \]
			\item[--] the set $\cM_{\rm inv}(g, {\Lambda_\e})$ is  contained  in the $\kappa\cdot(|\chi^c(\nu)|+\e)$-neighborhood of $\nu. $ 
		\end{itemize}
	Furthermore,  the set $\big\{\nu\in\mathcal{M}_{erg}(g)|~ \chi^c(\nu)\geq 0 \big\}$ is path connected.
	
  	\end{theoremalph}
  	\begin{Remark}
  		\begin{enumerate}
  			\item In Theorem~\ref{thm.approaching-non-hyperbolic-ergodic-measure}, if $\nu$ is non-hyperbolic (i.e. $\chi^c(\nu)=0$), then $\nu$ is approached by horseshoes in weak$*$-topology and in entropy;
  			\item  For any $g\in\cU$, all the periodic points with positive center Lyapunov exponent are homoclinically related\footnote{Recall that two hyperbolic periodic orbits are homoclinically related, if the stable manifold of one periodic orbit intersects the unstable manifold of the other transversely, and vice versa.} due to Theorem~\ref{thm.robust-minimal-ss}.
  			\item Comparing  with the results in \cite{BZ,DGS,YZ}, we  assume neither the   minimality of   strong unstable foliation    nor the existence of blenders.
  		\end{enumerate}
  		\end{Remark}

  		Theorem~\ref{thm.approaching-non-hyperbolic-ergodic-measure} is not only used to prove  Theorem~\ref{thm.existence} and has its own interest. The approximation of ergodic measures by hyperbolic sets   in various way has been studied. See for instance the results in $C^{1+\alpha}$-setting  by \cite{Ka,SW, Ge} and in $C^1$-setting with domination by \cite{G2,C,Ge} for hyperbolic ergodic measures, and the results in $C^1$-generic setting without domination by \cite{BCF} for ergodic measures.
  		The approximation of non-hyperbolic ergodic measures came into sight in recent years, and many results have been obtained. See for instance ~\cite{DGR1,BZ,YZ,DGS}, and one can refer to \cite{DG} for more results.
  		The main novelty here is that we neither assume the existence of blenders nor   the minimality of strong unstable foliation.  We use mostly expanding property to replace the existence of blenders and minimality of strong unstable foliation, and in some sense, mostly expanding property which is a non-uniformly hyperbolic property can play the same role in these problems.

  		Besides, one can get more information on the set of points  with vanishing center Lyapunov exponent, which in general is a larger set than the union of the generic points of non-hyperbolic ergodic measures. 
  		Let $f$ be a partially hyperbolic diffeomorphism with one-dimensional center bundle, and let us define the $0$-level  of  the center Lyapunov regular set:
  		\[\mathcal{R}_f(0):=\big\{x\in M| \lim_{|n|\rightarrow\infty}\frac{1}{n}\log\|Df^n|_{E^c(x)}\|=0 \big\}. \]
  		\begin{theoremalph}\label{thm.skeleton}
  			Let $f\in\diff^1(M)$ be a   partially hyperbolic diffeomorphism with $\dim(E^c)=1$. 
  		Assume that 
  		\begin{itemize}
  					\item $\chi^c(\mu)>0$ for any $\mu\in G^u(f);$
  			\item  the strong stable foliation is minimal.
  		\end{itemize}  
  		Then there exists a $C^1$-neighborhood  $\cU\subset\diff^1(M)$ of $f$ such that for  any $g\in\cU$, and any $h\leq h_{top}(g, \mathcal{R}_g(0))$ and $\e>0$,   there exists a hyperbolic basic set $\Lambda_\e$ of $g$ whose center is uniformly expanding such that 
  			\begin{itemize}
  				\item[--] $h_{top}(g,\Lambda_{\e})>h-\e, $
  				\item[--] for each   $\mu\in\mathcal{M}_{erg}(g,\Lambda_\e)$, one has $0<\chi^c(\mu,g)<\e.$
  			\end{itemize}
  		\end{theoremalph}
  	\begin{Remark}
  	For any non-hyperbolic ergodic measure $\nu\in\mathcal{M}_{erg}(g)$, one has $\nu(\mathcal{R}_g(0))=1$. By Theorem 3 in \cite{B} and the monotonicity of entropy, one has $h_\nu(g)\leq h_{top}(g, \mathcal{R}_g(0))$. Thus, Theorem~\ref{thm.skeleton} also implies that non-hyperbolic ergodic measures are approached by horseshoes in entropy. 
  	\end{Remark}
\noindent{\it Acknowledgments.} The author would like to thank C. Bonatti, L. D\'iaz,  S. Crovisier, K. Gelfert and  R. Saghin for helpful comments.   The author benefits a lot from the discussions with L. D\'iaz,  S. Crovisier, K. Gelfert,  D. Yang and  J. Yang.	The author is partially supported by National Key R\&D Program of China (2022YFA1005801),  National Key R$\&$D Program of China (2021YFA1001900),  NSFC 12001027 and the Fundamental Research Funds for the Central Universities.
  		
  \section{Preliminaries}
  In this section, we collect the results and notions that involved in this paper.
  
  \subsection{Topological entropy and metric entropy}\label{s.definition-of-topo-metric-entropy}
  In this section,   let $f:X\to X$ be a homeomorphism on  a compact metric space $(X,\ud)$.
  
  Given $n\in\mathbb{N}$ and $\e>0$, let us recall that 
  \begin{itemize}
  	\item  a \emph{$(n,\e)$-Bowen ball} centered at a point $x\in X$ is defined as 
  	\[B_n(x,\e)=\big\{y\in X: \ud(f^{i}(x),f^i(y))<\e \textrm{~for any $0\leq i<n$} \big\}. \]
  	\item a subset $S\subset X$ is called a \emph{$(n,\e)$-separated set}, if for any two different points $x,y\in S$, there exists $j\in\{0,\cdots,n-1 \}$ such that 
  	$\ud(f^j(x),f^j(y))>\e.$
  \end{itemize}
   Given an  invariant compact subset $K\subset X$, let $s(n,\e,K)$ be the maximal cardinality of $(n,\e)$-separated sets contained in $K$. 
   Then \emph{the topological entropy of $f$ on $K$} is defined as 
   \[ h_{top}(f,K)=\lim_{\e\rightarrow 0}\limsup_{n\rightarrow+\infty}\frac{1}{n}  \log s(n,\e,K) .\]
   When $K=X$, we simply denote $h_{top}(f):=h_{top}(f,X).$
   
  In analogy with the definition of Hausdorff dimension,  Bowen~\cite{B} introduced topological entropy for non-compact sets via open covers and here we present an equivalent one via Bowen balls (see \cite{Pe}). Let $Y\subset X$. For any $\e>0$ and $h\in\mathbb{R}$, let us define 
 \[  {\rm m}_{\e,h}(Y)=\lim_{n\rightarrow+\infty}\inf\big\{\sum_{i\in\mathbb{N}}e^{-hn_i}|~ \textrm{$Y\subset \cup_{i\in\mathbb{N}}B_{n_i}(x_i,\e)$ and $n_i\geq n$} \big \}.\]
 Define 
 \[h_{top}(f, Y,\e)=\inf\big\{h|~ {\rm m}_{\e,h}(Y)=0 \big\}=\sup\big\{h|~ {\rm m}_{\e,h}(Y)=+\infty\big\}. \]
 Now,   the topological entropy of $f$ on $Y$ is defined as 
 \[h_{top}(f,Y)=\lim_{\e\rightarrow 0}h_{top}(f, Y,\e). \]
  \begin{Remark}
 	Bowen showed that if  $Y\subset X$ is   compact, then the dimension-like definition of topological entropy coincides with the canonical definition using separated sets (see \cite[Proposition 1]{B}).
 \end{Remark}
 
 One can also define the entropy for  the non-compact set $Y$ by counting the cardinality of separated sets.  Let $s(n,\e,Y)$ be the maximal cardinality of $(n,\e)$-separated sets in $Y$, then one can define the \emph{lower capacity entropy}  and \emph{upper capacity entropy} of $f$ on $Y$ as follows:
 \[\underline{Ch}_{top}(f,Y)=\lim_{\e\rightarrow 0}\liminf_{n\rightarrow+\infty}\frac{1}{n}\log{s(n,\e,Y)}\]
 and 
 \[\overline{Ch}_{top}(f,Y)=\lim_{\e\rightarrow 0}\limsup_{n\rightarrow+\infty}\frac{1}{n}\log{s(n,\e,Y)}\]

 Now, we recall some basic properties for topological entropy.
 \begin{Proposition}[Proposition 2 in \cite{B} and Section 11 in \cite{Pe}]~\label{p.basic-property-entropy}
 	Let $f\in\homeo(X)$ be a homeomorphism on a compact metric space $(X,\ud)$. Then one has 
 	\begin{itemize}
 		\item for any subsets $Y_1\subset Y_2\subset X$, one has $h_{top}(f,Y_1)\leq h_{top}(f,Y_2).$
 		\item $h_{top}(f,Y)=\sup_n h_{top}(f,Y_n)$, where $Y=\cup_{n\in\mathbb{N}}Y_n\subset X.$
 		\item $\overline{Ch}_{top}(f,Y)\geq \underline{Ch}_{top}(f,Y)\geq h_{top}(f,Y).$
 	\end{itemize} 
 \end{Proposition}

   \medskip

   At last, we recall the definition of metric entropy of  an $f$-ergodic measure $\mu$ given by Katok~\cite{Ka}. 
   Given $\delta\in(0,1)$, $n\in\mathbb{N}$ and $\e>0$, let $s(n,\e,\delta)$ be the minimal cardinality of $(n,\e)$-Bowen balls whose union covers a set with $\mu$-measure no less than $\delta.$ Then the \emph{metric entropy of $f$ with respect to $\mu$} is defined as 
   \[h_\mu(f)=\lim_{\e\rightarrow 0}\limsup_{n\rightarrow+\infty}\frac{1}{n}\log s(n,\e,\delta) =\lim_{\e\rightarrow 0}\liminf_{n\rightarrow+\infty}\frac{1}{n}\log s(n,\e,\delta).\]
   \subsection{Liao's shadowing lemma and Pliss lemma}
   A $Df$-invariant splitting $T_\La M=E\oplus F$ over an invariant  compact   set $\Lambda$ of $f\in\diff^1(M)$ is a \emph{dominated splitting}, if there exists $N\in\mathbb{N}$ such that 
   \[\|Df^N|_{E(x)}\|\cdot \|Df^{-N}|_{F(f^N(x))}\|\leq \frac 12 \textrm{~for any $x\in\Lambda$.} \]
   
  Now, we recall Liao's shadowing lemma  which was improved  by S. Gan.  This shadowing lemma allows us to find periodic orbits chasing (long) orbit segments with "weak" center Lyapunov exponent for a large proportion of time. 

 \begin{Theorem}[\cite{L,G}]~\label{l.liao-gan-shadowing}
	Let $f\in\diff^1(M)$ and $\Lambda$ be an $f$-invariant compact set. Assume that 
	$\Lambda$ admits a dominated splitting of the form $T_\Lambda M=E\oplus F$.
	Then  for any $\lambda\in(0,1)$, there exist constants $L>1$ and $d_0>0$  such that for any $d\in(0,d_0)$, any $x\in\Lambda$ and any $n\in\mathbb{N}$  satisfying  
	\begin{itemize}
		\item $\ud(f^{n}(x),x)<d;$  
		\item \[\|Df^j|_{E(x)}\|<\lambda^j\textrm{~and~} \|Df^{-j}|_{F(f^{n}(x))}\|<\lambda^j, 
		\textrm{~for any  $1\leq j\leq n$};\]
	\end{itemize}
	there exists a periodic point $p$ of period $n$ such that 
	\[ \ud(f^i(x),f^i(p))<L\cdot d  \textrm{~ for any $0\leq i\leq n-1$}.\]
\end{Theorem} 
  \begin{Remark}
  	We will apply this shadowing lemma to the splitting $E^s\oplus(E^c\oplus E^u)$ and in this case, one only needs to consider   the norm along the bundle $E^c\oplus E^u.$
  	\end{Remark}
  Now, we recall the Pliss lemma  which can be used for finding points with uniform size of stable or unstable manifolds.
  \begin{Lemma}[Pliss lemma \cite{Pl}]\label{l.pliss}
  	Let $a_1,\cdots, a_k$ be some real numbers and assume that $\max_{i\leq k} a_i\leq C$ for some   $C\in\RR$. Suppose that $\sum_{i=1}^k a_i \geq k\chi_1$ for some $\chi_1$. Then for any $\chi_2<\chi_1$, there exist $1\leq j_1<j_2<\cdots<j_l\leq k$ such that 
  	\begin{itemize}
  		\item $\rho:=\frac{l}{k}\geq \frac{\chi_1-\chi_2}{C-\chi_2}$;
  		\item for each $1\leq n\leq l$ and each $1\leq m\leq j_n$, one has \[\sum_{i=m}^{j_n} a_i\geq (j_n-m+1)\chi_2.\]
  	\end{itemize}
  	
  \end{Lemma}
 
  \subsection{Plaque families and uniform size of invariant  manifolds}
  
   The existence of plaque family was given by Theorem 5.5 in \cite{HPS} for a single diffeomorphism  and was extended to a neighborhood of a diffeomorphism. See Lemma 3.5 in \cite{CP}. Given a diffeomorphism $f\in\diff^1(M)$, let $\La$ be an invariant  compact   set and $E$ be a vector bundle over $\Lambda$.     For $x\in\La$ and $r>0$, let us denote 
   \[E(x,r):=\big\{v\in E(x)|~ \|v\|<r \big\} \textrm{~and~} E(r)=\cup_{x\in\La} E(x,r).\]
   
   \begin{Theorem}[Plaque Family Theorem]\label{thm.plaque-family}
   Let  $f\in\diff^1(M)$  and  $\La$  be an invariant compact  set with the dominated splitting of the form  $T_\La M=E\oplus F.$     
   
   Then there exist a $C^1$-neighborhood $\cU$ of $f$ and a neighborhood $U$ of $\Lambda$ such that for any $g\in\cU$,  the maximal invariant compact  set $\La_g$ of $g$ in $U$ admits a dominated splitting $T_{\La_g}M=E_g\oplus F_g,$  and there exist two continuous families of  maps $\cW^{cs}_g: E_g(1)\to M$ and $\cW^{cu}_g: F_g(1)\to M$ satisfying the following properties:
   \begin{itemize}
   	\item  for each $x\in\La_g$, the map $\cW^{cs}_{x,g}:E_g(x,1)\to M$ (resp. $\cW^{cu}_{x,g}:F_g(x,1)\to M$) is a $C^1$-embedding, $x=\cW^{cs}_{x,g}(0_x)$ (resp. $x=\cW^{cu}_{x,g}(0_x)$ ) and its  graph  is tangent to $E_g(x)$ (resp. $F_g(x)$) at the point $x$;
 \item   	the families $\{\cW^{cs}_{x,g} \}$ and $\{\cW^{cu}_{x,g} \}$ of $C^1$-embedding maps are  continuous with respect to $x,g$ in $C^1$-topology; 
 \item for any $\delta\in(0,1)$, there exists $\delta'>0$ such that 
 \[g(\cW^{cs}_{x,g}(E_g(x,\delta')))\subset \cW^{cs}_{g(x),g}(E_g(g(x),\delta))\] 
 and  
 \[  g^{-1}(\cW^{cu}_{x,g}(F_g(x,\delta')))\subset \cW^{cu}_{g^{-1}(x),g}(F_g(g^{-1}(x),\delta)).\]
   \end{itemize}
   For $g\in\cU$, $x\in\La_g$ and $\delta\in(0,1)$, we denote  by $\cW^{cs}_\delta(x,g):=\cW^{cs}_{x,g}(E_g(x,\delta)))$ and 
   $\cW^{cu}_\delta(x,g):=\cW^{cu}_{x,g}(F_g(x,\delta)))$.
   	\end{Theorem}
   	When there is no ambiguity, we will drop the index $g$ for simplicity. The $C^1$-submanifolds $\cW^{cs}_\delta(x,g)$ (resp. $\cW^{cu}_\delta(x,g)$) is called the $cs$-plaque (resp.  $cu$-plaque) centered at $x$ and of radius $\delta$. The last item in Theorem~\ref{thm.plaque-family} tells us that the $cu$-plaques and $cs$-plaques are locally $g$-invariant.

   The following result gives the uniformity on the size of   unstable manifolds for a single diffeomorphism, whose proof can be found in \cite[Section 8.2]{ABC}.
   \begin{Lemma}\label{l.uniform-size-of-unstable-for-single-diffeo} 
   	Let  $f\in\diff^1(M)$  and  $\La$  be an invariant  compact   set  with the dominated splitting of the form  $T_\La M=E\oplus F.$  Consider a $cu$-plaque family $\cW^{cu}$ corresponding to the bundle $F$. 
   	For any $\chi>0$, there exists $\delta>0$ such that  if   $x\in \La$ satisfies that 
   	\[\prod_{i=0}^{n-1}\|Df|_{F(f^{-i}(x))}\|\leq e^{-n\chi} \textrm{~for any $n\in\mathbb{N}$}, \]
   	then $\cW^{cu}_{\delta}(x)$ is contained in the unstable manifold of $x$.
   \end{Lemma}

  Using the uniform continuity of the plaque families among nearby systems   and Lemma~\ref{l.uniform-size-of-unstable-for-single-diffeo}, one can obtain the similar result for diffeomorphisms close to $f$ and for our purpose we state it in partially hyperbolic setting. 
  \begin{Lemma}~\label{l.uniform-unstale-manifold}
  	Let $f\in\diff^1(M)$ be a partially hyperbolic diffeomorphism with   $\dim(E^c)=1$, and $\widetilde{\cU}$ be a $C^1$-neighborhood of $f$ given by Lemma~\ref{thm.plaque-family} together with a plaque family  $\cW^{cu}_g$ with $g\in\widetilde\cU$,  which corresponds  to the bundle $E_g^c\oplus E_g^u$ and depends continuously on $g$. 
  	
  	Then for any $\chi>0$, there exist  a $C^1$-small neighborhood $\cU\subset \widetilde{\cU}$ of $f$,  and a small constant $\e_0>0$ such that 
  	\begin{itemize}
  		\item for  any $x\in M$, one has
  		\begin{itemize}
  			\item[--] $g^{-1}(\cW^{cu}_{\e_0}(x,g))\subset \cW^{cu}_{1/2}(g^{-1}(x),g)$;
  			\item[--] $\big|\log\|Dg^{-1}|_{E_g^c(x)}\|-\log\|Dg^{-1}|_{T_y\cW^{cu}_{\e_0}(x,g)}\|\big|<\chi/4\, \textrm{~for any $y\in \cW^{cu}_{\e_0}(x,g)$};$
  		\end{itemize} 
  		\item 	 if a point $x\in M$ satisfies that $\|Dg^{-n}|_{E_g^c(x)}\|<e^{-n\chi}$ for any $n\in\NN$, then   
  		\begin{itemize}
  			\item $\cW^{cu}_{\e_0}(x,g)$ is contained in the unstable manifold of $x$ and is tangent everywhere to $E_g^c\oplus E_g^u;$
  			\item   $\|Dg^{-n}|_{E_g^c(y)}\|<e^{-3n\chi/4}$ for any $n\in\NN$ and any $y\in \cW^{cu}_{\e_0}(x,g)$.
  		\end{itemize}
  	\end{itemize}
  \end{Lemma}

   Given $\e>0$, $\ell>0$ and a plaque family $\cW^{cu}$ corresponding to the bundle $E^c\oplus E^u$, one says that the strong stable foliation $\mathcal{F}^s$ is \emph{$(\ell,\e)$-dense with respect to   $\cW^{cu}$}, 
  if for any $x,y\in M$ the local strong stable manifold $\mathcal{F}^s_{\ell}(x)$ has non-empty transverse intersection with $\cW^{cu}_\e(y)$, where $\mathcal{F}^s_{\ell}(x)$ is the $\ell$-neighborhood  of $x$ in $\mathcal{F}^s(x)$ under its intrinsic topology.
  
  The following result comes from the uniform continuity of strong stable manifolds and the plaque families with respect to the bundle $E^c\oplus E^u$.
  \begin{Lemma}\label{l.minimal-strong-stable-to-cu-dense}
  	Let $f\in\diff^1(M)$ be a partially hyperbolic diffeomorphism and $\widetilde{\cU}$ be a $C^1$ neighborhood of $f$ given by Lemma~\ref{thm.plaque-family} together with a plaque family  $\cW^{cu}_g$ with $g\in\widetilde\cU$,  which corresponds  to the bundle $E_g^c\oplus E_g^u$ and depends continuously on $g$. 
  	 Assume, in addition,  that   strong stable foliation of $f$ is minimal.

  	  Then for any $\e>0$, there exist  a $C^1$-neighborhood $\cU_\e\subset\widetilde\cU$ of $f$ and a constant  $\ell>0$ such that the strong stable foliation of $g\in\cU_\e$ is $(\ell,\e)$-dense with respect to  $\cW_g^{cu}$. 
  \end{Lemma} 
\subsection{Constructions of horseshoes from hyperbolic periodic orbits}
  Katok and Mendoza  \cite{KM} gave a way to construct horseshoes approaching  hyperbolic ergodic measures on surfaces   in $C^{1+\alpha}$-setting and this has been generalized to higher dimension in \cite{Ge}. Here we state a mechanism using a collection of periodic orbits to find horseshoes  in  $C^1$-setting with domination. 
  
 Given an invariant  compact set $\Lambda$ of $f\in\diff^1(M)$ exhibiting    a dominated splitting of the form $T_\Lambda M=E\oplus F$  and given $\chi>0$, let us define 
  \[{\rm NUH}_\chi=\big\{x\in\Lambda|\, \prod_{i=0}^{k-1}\|Df|_{E(f^i(x))}\|\leq e^{-k\chi},~\prod_{i=0}^{k-1}\|Df^{-1}|_{F(f^{-i}(x))}\|\leq e^{-k\chi} \textrm{~for any $k\in\NN$}  \big\}. \]
  
  \begin{Theorem}[Theorem 4.1 in \cite{YZ}]\label{thm.existence-of-horseshoes}
  	Let $f\in\diff^1(M)$ and $\Lambda$ be an invariant compact set admitting a dominated splitting of the form $T_\Lambda M=E\oplus F$. For any $\chi>0$
  	and  any $\e>0$, there exists $\xi_0>0$ satisfying the following properties. For any $n\in\NN$ and any $\xi\in(0,\xi_0)$, if there exist periodic points $p_1,\cdots,p_m\in{\rm NUH}_\chi$ of the same period $l$ such that 
  	\begin{itemize}
  		\item $d(p_i,p_j)<\xi/16$ for $1\leq i<j\leq m$;
  		\item $\{p_1,\cdots,p_m \}$ is a $(l, \xi)$-separated set;
  	\end{itemize}
  	then there exists a   hyperbolic basic set $K$ whose stable bundle has dimension $\dim(E)$ such that 
  	\begin{itemize}
  		\item[--]  $h_{top}(f,K)\geq \frac{\log m}{l}$;
  		\item[--]   the set $\cM_{\rm inv}(f, K)$ is contained in the $\e$-neighborhood of $\big\{\sum_{i=1}^mt_i\delta_{\cO_{p_i}}|
  		\,t_i\geq 0, ~ \sum_{i=1}^m t_i=1 \big\}$, where $\delta_{\cO_{p_i}}=\frac{1}{l}\sum_{j=0}^{l-1}\delta_{f^j(p_i)}.$ 
  	\end{itemize}
  \end{Theorem}

  \subsection{Unstable entropy and entropy formulas}
 In this section, we firstly recall the definition of  the metric entropy along unstable  manifolds, which    was  defined by F. Ledrappier and L. Young \cite{LeYo1,LeYo2}. Then we recall its applications based on some entropy formulas. 
  
  
  Given a Borel probability measure $\mu$ and a foliation $\cF$ on the manifold $M$, a measurable partition $\cA$ is called \emph{$\mu$-subordinate to $\cF$} if for $\mu$ a.e. $x\in M$, one has 
  \begin{itemize}
  	\item $\cA(x)\subset \cF(x)$, where $\cA(x)$ denotes the element of $\cA$ containing $x$;
  	\item  there exists $r_x>0$ such that $\cF_{r_x}(x)\subset \cA(x)$, where $\cF_{r_x}(x)$ denotes the $r_x$-neighborhood of $x$ in $\cF(x)$ under its intrinsic topology.
  \end{itemize}
  
  Given $f\in\diff^1(M)$, a measurable partition $\cA$ is \emph{increasing} if $f(\cA)\prec\cA$. If $f\in\diff^1(M)$ is partially hyperbolic and $\mu$ is an $f$-invariant measures, then    there exists an increasing measurable partition which is $\mu$-subordinate to the strong unstable foliation due to \cite{LeS} (see also \cite{Y}). Given two increasing measurable partitions $\cA_1,\cA_2$ which are $\mu$-subordinate to the unstable foliation, Ledrappier and Young \cite{LeYo1} proved that $H_\mu(\cA_1|f(\cA_1))=H_\mu(\cA_2|f(\cA_2))$, and  this yields the following definition.
  \begin{Definition}
  	Let $f\in\diff^1(M)$ be a partially hyperbolic diffeomorphism and $\mu\in\cM_{inv}(f)$. The \emph{unstable (metric)-entropy} of $\mu$ is defined as 
  	\[h_\mu(f,\cF^u)=H_\mu\big(\cA|f(\cA)\big), \]
  where $\cA$	is an increasing measurable partition  $\mu$-subordinate to $\cF^u.$
  \end{Definition}

  Recall that an invariant measure of a $C^{1+\alpha}$-partially hyperbolic diffeomorphism is called a \emph{$u$-Gibbs state}, if its conditional measures along the strong unstable manifolds are absolutely continuous with respect to the Lebesgue measure on the strong unstable manifolds. The existence of $u$-Gibbs states was firstly proved by Y. Pesin and Y. Sinai \cite{PS}. It has been shown by Ledrappier that   u-Gibbs property is equivalent to the unstable  entropy formula.  
  \begin{Theorem}[Th\'eor\`eme in \cite{Le}]\label{thm.ledrappier}
 Let $f\in\diff^{1+\alpha}(M)$ be partially hyperbolic and $\mu\in\cM_{inv}(f)$. 
 Then $\mu$ is a $u$-Gibbs state if and only if $h_\mu(f,\mathcal{F}^u)=\int\log|\det(Df|_{E^u})|\ud\mu.$
  \end{Theorem}

For a $C^1$-partially hyperbolic diffeomorphism $f$, let us denote 
\begin{equation}\label{eq:unstable-entropy-formula}
G^u(f):=\big\{\mu\in\mathcal{M}_{inv}(f)|h_\mu(f,\mathcal{F}^u)=\int\log|\det(Df|_{E^u})|\ud\mu\big\}.
\end{equation}
The following result describes the set $G^u(f)$.
\begin{Lemma}[\cite{CYZ,HYY,Y}]\label{l.the-set-Gu(f)}
 Let $f$ be a $C^1$-partially hyperbolic diffeomorphism. Then
\begin{itemize}
	\item $G^u(f)$ is a non-empty convex compact set;
	\item $G^u(f)$ varies upper semi-continuously with respect to $f$;
	\item if $\mu$ belongs to $G^u(f)$, so do its ergodic components.
\end{itemize}
\end{Lemma}
In $C^1$-setting, the large deviation property for the set $G^u(f)$ holds due to \cite{CYZ}. Recall  that for the partially hyperbolic splitting $TM=E^s\oplus E^c\oplus E^u$, one can associate a continuous cone field of size $\beta>0$ to  the strong unstable bundle $E^u$: 
\[\mathcal{C}^u_\beta:=\big\{v\in TM| ~v=v^{cs}+v^u,~ v^{cs}\in E^s\oplus E^c, ~v^u\in E^u,~ \|v^{cs}\|\leq \beta \cdot\|v^{u}\|\big\} .\]
For $\beta>0$ small, the cone field $\mathcal{C}^u_{\beta}$ is strictly $Df$-invariant, that is, there exists $\beta'\in(0,\beta)$ such that $Df(\mathcal{C}^u_\beta(x))\subset \mathcal{C}^u_{\beta'}(f(x))$ for any $x\in M$.
When there is no ambiguity, we would drop the index $\beta$ for simplicity. 
\begin{Theorem}[Theorem D$^\prime$ in \cite{CYZ}]\label{thm.large-deviation}
	Let $f\in\diff^1(M)$ be a partially hyperbolic diffeomorphism. Then for any continuous function $\varphi:M\to\mathbb{R}$ and any $\kappa>0$, there exist positive constants  $r_0,a_\kappa,b_\kappa>0$ and a $Df$-invariant cone field $\mathcal{C}^u$ for the strong unstable bundle such that for any  disc $D$ tangent to the cone field $\mathcal{C}^u$ and of diameter smaller than $r_0$, one has 
	\[\Leb_D\bigg(\big\{x\in D: ~~\ud\big(\frac{1}{n}S_n\varphi(x), I(\varphi)\big)>\kappa\big\}\bigg)<a_\kappa\cdot e^{-n b_\kappa} \textrm {~ for any $n\in\mathbb{N}$},\]
	where $\Leb_D$ denotes the Lebesgue measure on the submanifold $D$, $S_n\varphi=\sum_{i=0}^{n-1}\varphi\circ f^i$ and \[I(\varphi)=\big\{\int\varphi\ud\nu: \nu\in G^u(f) \big\}.\]
\end{Theorem}

In $C^{1+\alpha}$-setting, $u$-Gibbs states are the candidate for SRB measures. Recall that an invariant measure is called an \emph{SRB} if it admits positive Lyapunov exponents and its conditional measures along the Pesin's unstable manifolds are absolutely continuous with respect to the Lebesgue measure. 
\begin{Theorem}[Theorem A in \cite{LeYo1}]\label{thm.Ledrappier-Young}
Let $f\in \diff^{1+\alpha}(M)$ and $\mu$ be an invariant measure with positive Lyapunov exponents. Then $\mu$ is an SRB if and only if $\mu$ satisfies Pesin's entropy formula, that is,
\[h_\mu(f)=\int\sum\lambda^+(x)\ud\mu, \]
where $\sum\lambda^+(x)$ is the sum of positive Lyapunov exponents of $x$ with multiplicities counted.
\end{Theorem}
\begin{Remark}
	Ledrappier and Young \cite{LeYo1} proved this result in $C^2$-setting and Brown \cite{Br} generalized it to $C^{1+\alpha}$-setting.
	\end{Remark}
In \cite{LeYo2}, Ledrappier and Young also gave an entropy formula in terms of Lyapunov exponents and transverse dimensions. The following result is a direct consequence of Corollary 7.2.2 in \cite{LeYo2} which assumes that $f$ is $C^2$, and it has been relaxed to $C^{1+\alpha}$-setting in \cite{Br}. 
\begin{Theorem}\label{thm.ledrappier-young2}
	Let $f\in\diff^{1+\alpha}(M)$ be a partially hyperbolic diffeomorphism and $\mu\in\mathcal{M}_{inv}(f)$. If all the center Lyapunov exponents of $\mu$ are non-positive, then 
	$h_\mu(f)=h_\mu(f,\cF^u). $
\end{Theorem}
 \subsection{Mostly expanding partially hyperbolic diffeomorphisms}~\label{s.mostly-expanding}
Recall  that a partially hyperbolic diffeomorphism $f\in\diff^{1+\alpha}(M)$ is \emph{mostly expanding}, if all the center Lyapunov exponents of each $u$-Gibbs state  of $f$ are positive. 

The following result comes from Lemma 4.1 and Theorem B in \cite{AV}, and one can also refer to the proof of  Theorem C in \cite{Y}. It tells us that mostly expanding is a $C^1$-open property.
\begin{Theorem}\label{thm.robust-mostly-expanding}
	Let $f\in\diff^{1+\alpha}(M)$ be a mostly expanding partially hyperbolic diffeomorphism. Then there exist a constant $\chi>0$, an integer $n_0\in\mathbb{N}$ and  a $C^1$-neighborhood $\cV$ of $f$ such that for each $g\in\cV\cap\diff^{1+\alpha}(M)$ and any $u$-Gibbs state $\mu$ of $g$, one has 
	\[\int\log\|Dg^{-n_0}|_{E^c}\|\ud\mu<-\chi. \]
\end{Theorem}

As the set of $u$-Gibbs states of $f\in\diff^{1+\alpha}(M)$ coincides  with the compact set $G^u(f)$  (due to Theorem~\ref{thm.ledrappier})  and  
$G^u(f)$ varies upper semi-continuously with respect to $f$ (due to Lemma~\ref{l.the-set-Gu(f)}), one immediately obtains the following.
\begin{Corollary}\label{cor:robust-mostly-expanding}
	Let $f\in\diff^{1+\alpha}(M)$ be a mostly expanding partially hyperbolic diffeomorphism. Then there exist a constant $\chi>0$, an integer $n_0\in\mathbb{N}$ and  a $C^1$-neighborhood $\cV$ of $f$ such that for each $g\in\cV$ and any   $\mu\in G^u(g)$, 
	one has 
	\[\int\log\|Dg^{-{n_0}}|_{E^c}\|\ud\mu<-\chi.\]
\end{Corollary}

In general the minimality of strong stable foliation is not an open property. Nevertheless, the following result tells us that under the hypothesis of mostly expanding property,  the minimality of strong stable foliation is a $C^1$-open property.  
\begin{Theorem}\label{thm.robust-minimal-ss}
	Let $f\in\diff^1(M)$ be a partially hyperbolic diffeomorphism with $\dim(E^c)=1.$ 
	Assume that 
	\begin{itemize}
		\item $\int\log\|Df|_{E^c}\|\ud\mu>0$ for any $\mu\in G^u(f)$;
		\item the strong stable foliation is minimal.
	\end{itemize}
	Then there exists a $C^1$-neighborhood $\mathcal{U}\subset\diff^1(M)$ of $f$ such that for each $g\in\cU$, one has
	\begin{itemize}
		\item[--] $\int\log\|Dg|_{E_g^c}\|\ud\mu>0$ for any $\mu\in G^u(g)$;
		\item[--]  the strong stable foliation of $g$ is  minimal.
		\end{itemize}
\end{Theorem}
\proof
As $\dim(E^c)=1$, by Corollary~\ref{cor:robust-mostly-expanding},
there exist  a $C^1$-neighborhood $\cV_1$ of $f$ and a constant $\chi>0$ such that for any $g\in\cV_1$, each invariant measure $\mu\in G^u(g)$ satisfies $\int\log\|Dg|_{E^c_g}\|\ud\mu>\chi.$ For the constant $\chi>0$, by Lemma~\ref{l.uniform-unstale-manifold}, there exist a $C^1$-neighborhood $\widetilde\cV_2\subset\diff^1(M)$ of $f$, a plaque family $\cW^{cu}_g$ with $g\in\widetilde\cV_2$ and a constant $\e_0>0$ satisfying the posited properties. 
As the strong stable foliation of $f$ is minimal, by Lemma~\ref{l.minimal-strong-stable-to-cu-dense},  there exist a $C^1$-neighborhood $\cV_2\subset\widetilde\cV_2$ of $f$ and $\ell_0>0$ such that 
for any $g\in\cV_2$ and any $x,y\in M$, the local strong stable manifold $\mathcal{F}^s_{\ell_0}(x)$ of $g$ has  non-empty transverse intersection with $\cW_{\e_0}^{cu}(y,g).$

Let $\cU=\cV_1\cap\cV_2$, and fix $g\in\cU$. Let us denote by  $\chi^u:=\inf_{\mu\in G^u(g)}\int\log\|Dg|_{E^c_g}\|\ud\mu>\chi,$  $\varphi(x)=\log\|Dg|_{E^c_g(x)}\|$ and $\kappa=\frac{\chi^u-\chi}{2}.$ Now, apply Theorem~\ref{thm.large-deviation} to the continuous function  $\varphi$ and $\kappa$, and one obtains positive constants $r_0, a_\kappa, b_\kappa$ such that for any $w\in M$, one has 
\begin{equation}\label{eq:large-deviation-ss-minimal}
	\Leb_{\mathcal{F}_{r_0}^u(w)}\bigg(\big\{y\in \mathcal{F}_{r_0}^u(w): ~~\ud\big(\frac{1}{n}S_n\varphi(y), I(\varphi)\big)>\kappa\big\}\bigg)<a_\kappa\cdot e^{-n b_\kappa}\, \textrm{~for any $n\in\mathbb{N}.$} 
	\end{equation}

Let $\chi_{\rm max}(g)=\sup_{x\in M}|\varphi(x)|.$ Apply $C=\chi_{\rm max}(g)$, $\chi_1=(\chi^u+\chi)/2$ and $\chi_2=\chi$ to Lemma~\ref{l.pliss}, and one obtains $\rho\in(0,1)$ with the posited properties

Fix $x,z\in M$. We are going to show that $\mathcal{F}^s(z)\cap B_{2\e}(x)\neq\emptyset$ for any $\e>0$. 
For  any $\e>0$, by the uniform expansion of $g$ along the strong unstable foliation, there exists an integer $k=k(\e)$ such that $\mathcal{F}^u_{r_0}(g^k(x))\subset g^k(\mathcal{F}_\e^u(x)).$ Take $N_1$ large enough such that for any $j\geq N_1,$  one has
\begin{itemize}
	\item  $a_\kappa\cdot e^{-j b_\kappa}< \Leb_{\mathcal{F}_{r_0}^u(g^k(x))}(\mathcal{F}_{r_0}^u(g^k(x))/2,$ 
	\item  $e^{-j\chi/2}\cdot \e_0\cdot \|Dg^{-1}\|^k<\e$, where  $\|Dg^{-1}\|=\sup_{w\in M}\|Dg^{-1}(w)\|$.
\end{itemize}
For $n>\rho^{-1}N_1+1$, by Equation~\eqref{eq:large-deviation-ss-minimal}, one has 
\[	\Leb_{\mathcal{F}_{r_0}^u(g^k(x))}\bigg(\big\{y\in \mathcal{F}_{r_0}^u(g^k(x)): ~~\ud\big(\frac{1}{n}S_n\varphi(y), I(\varphi)\big)>\kappa\big\}\bigg)< \Leb_{\mathcal{F}_{r_0}^u(g^k(x))}(\mathcal{F}_{r_0}^u(g^k(x))/2.\]
Therefore there exists a point $y\in \mathcal{F}^u_{r_0}(g^k(x))$ such that $\ud\big(\frac{1}{n}S_n\varphi(y), I(\varphi)\big)\leq \kappa$, and thus $\frac{1}{n}S_n\varphi(y)\geq \chi^u-\kappa=\frac{\chi^u+\chi}{2}$. By Pliss lemma,  there exists $n'\geq n\rho> N_1$ such that 
\begin{equation}\label{eq:pliss-at-y-robut-s-minimal}
	\|Dg^{-j}|_{E^c_g(g^{n'}(y))}\|\leq e^{-j\chi} \, \textrm{~for any $j\in[1,n'].$}\end{equation}
As  $g\in\cU\subset\cV_2$, the local stable manifold $\mathcal{F}^s_ {\ell_0}(g^{k+n'}(z))$ has non-empty transverse intersection with $\cW^{cu}_{\e_0}(g^{n'}(y),g).$ As the strong stable foliation is $g$-invariant,   $\mathcal{F}^s(g^{k}(z))$ has non-empty transverse intersection with $g^{-n'}\big(\cW^{cu}_{\e_0}(g^{n'}(y),g)\big)$. As $g\in\cV_2\subset\widetilde\cV_2$, by Equation~\eqref{eq:pliss-at-y-robut-s-minimal}, the first item in  Lemma~\ref{l.uniform-unstale-manifold}  and the locally  invariant property of the plaque family,  one has \[g^{-n'}\big(\cW^{cu}_{\e_0}(g^{n'}(y),g)\big)\subset \cW^{cu}_{e^{-n'\chi/2}\e_0}(y,g).\]
 Note that $e^{-n'\chi/2}\cdot \e_0\cdot \|Dg^{-1}\|^{k}<\e$ due to $n'>N_1$, then   the strong stable manifold $\mathcal{F}^s(z)$ intersects   $\cW^{cu}_{\e}(g^{-k}(y),g)$ which is contained in $B_\e(g^{-k}(y))$. Since $g^{-k}(y)\in\cF^u_\e(x)$, one has $\mathcal{F}^s(z)\cap B_{2\e}(x)\neq\emptyset.$ By the arbitrariness of $x,z$ and $\e$,  one deduces that the strong stable foliation  of $g$ is minimal.
\endproof
It has been shown in \cite{Y} that there are plenty of mostly expanding partially hyperbolic diffeomorphisms among conservative  partially hyperbolic diffeomorphisms with positive center Lyapunov exponents. Recall that a partially hyperbolic diffeomorphism is \emph{accessible}, if any pair of points can be connected by a path which is a concatenation of paths lying entirely in a strong stable or  a strong unstable manifold. 
\begin{Theorem}[Theorem D in \cite{Y}]
	Let $f\in\diff_{\rm m}^{1+\alpha}(M)$ be a volume preserving partially hyperbolic diffeomorphism with $\dim(E^c)=1.$ 
	Assume that $f$ is accessible and the volume measure has positive center Lyapunov exponent. Then $f$ is mostly expanding. 
\end{Theorem}
 Instead of assuming accessibility, one can assume  the minimality of strong stable foliation and obtain the similar result. It will be used to prove Corollary~\ref{cor:conservative-case}.
 \begin{Theorem}\label{thm.citerion-mostly-expanding}
 	Let $f\in\diff_{\rm m}^{1+\alpha}(M)$ be a   partially hyperbolic diffeomorphism preserving a smooth volume ${\rm m}$.   Assume that
 	\begin{itemize}
 			\item   the strong stable foliation of $f$ is minimal;
 				\item    $\dim(E^c)=1$   and $\int\log\|Df|_{E^c}\|\ud \m>0$.
 	\end{itemize}  
 	Then one has that 
 	\begin{itemize}
 		\item[--] $f$ is mostly expanding;
 		\item[--]$(f, \m)$ is $C^1$-stably ergodic, that is, there exists a $C^1$-neighborhood $\mathcal{U}\subset\diff^1(M)$ of $f$  such that each $g\in\cU\cap \diff_{\rm m}^{1+\alpha}(M)$ is ergodic with respect to the volume measure $\m.$
 	\end{itemize}
 \end{Theorem}
 \proof 
 By the absolute continuity of the strong stable and unstable foliations, each ergodic component of the volume measure $\m$ is both a $u$-Gibbs state and an $s$-Gibbs state (i.e. a $u$-Gibbs state for $f^{-1}$). 
 
 By assumption, there exists an ergodic component $\m_0$ of $\m$ such that $\chi^c(\m_0,f)>0.$ Note that  $\m_0$ is a $u$-Gibbs state for $f$ as well as a $u$-Gibbs state for $f^{-1}$. As $\m_0$ has negative center Lyapunov exponent for $f^{-1}$, thus $\m_0$ is also an SRB for $f^{-1}$.
 \begin{Claim-numbered}\label{cl.unique-u-Gibbs}
 	$f^{-1}$ has a unique $u$-Gibbs state $\m_0$.
 	\end{Claim-numbered}
 \proof 
Let  $\mu$ be an ergodic $u$-Gibbs state for $f^{-1}$. Note that the conditional measures of $\m_0$ along strong unstable manifolds of $f^{-1}$ are absolutely continuous with respect to the Lebesgue measure, and $\m_0$ has negative center Lyapunov exponent for $f^{-1}$,  by the absolute continuity of the Pesin stable manifold for $(\m_0,f^{-1})$,   the minimality of the strong unstable foliation for $f^{-1}$ and the ergodicity of the measures $\m_0$ and $\mu$, one deduces that $\m_0=\mu$.
 \endproof
Claim~\ref{cl.unique-u-Gibbs} implies that $\m$ has only one ergodic component $\m_0$ proving the ergodicity of $\m.$
 

 Consider an ergodic $u$-Gibbs state $\mu$ for $f$,  and we need to show that $\mu$ has positive center Lyapunov exponent for $f$. By Theorem~\ref{thm.ledrappier}, one has $h_\mu(f,\mathcal{F}^u)=\int\log\det(Df|_{E^u})\ud\mu$. 
 If $\chi^c(\mu,f)\leq 0$, by Theorem~\ref{thm.ledrappier-young2}, one has 
 \[h_\mu(f)=h_\mu(f,\mathcal{F}^u)=\int\log\det(Df|_{E^u})\ud\mu.\]
 As $f$ is volume preserving, one has $\int\log\det(Df|_{E^u})\ud\mu=-\int\log\det(Df|_{E^s\oplus E^c})\ud\mu.$ This gives that   $h_\mu(f)=-\int\log\det(Df|_{E^s\oplus E^c})\ud\mu$. By Theorem~\ref{thm.Ledrappier-Young}, $\mu$ is an SRB measure for $f^{-1}$. By the absolute continuity of the strong unstable manifolds of $f^{-1}$, $\mu$ is also a $u$-Gibbs state for $f^{-1}$. Then Claim~\ref{cl.unique-u-Gibbs} gives that $\mu=\m_0$ which contradicts the fact that   $\chi^c(\mu,f)\leq 0$ and $\chi^c(\m_0,f)>0$. This proves that all the ergodic $u$-Gibbs states of $f$ have positive center Lyapunov exponent, and hence $f$ is mostly expanding.

Let $\cU\subset \diff^1(M)$ be the $C^1$-open neighborhood of $f$ given by Theorem~\ref{thm.robust-minimal-ss}. For each $g\in \cU\cap\diff^{1+\alpha}_{\m}(M),$ by the absolute continuity of the strong unstable foliation, one has   $\m\in G^u(g)$, and thus   $\int\log\|Dg|_{E^c_g}\|\ud\m>0$ due to the first item of Theorem~\ref{thm.robust-minimal-ss}.  Now apply the arguments above to $(g,\m)$ with $g\in\cU\cap\diff^{1+\alpha}_{\m}(M)$, and one deduces that $(g,\m)$ is ergodic. 
 \endproof
 \section{Existence  of periodic orbits  with prescribed behavior} 
This section gives the main ingredient for proving our main results. 
 We will first   find some periodic orbits which shadow the generic points of an  ergodic measure with non-positive center Lyapunov exponent for some proportion of time. 
 
Throughout this section, given a foliation $\mathcal{F}$ on $M$ with $C^1$-leaves, $x\in M$ and $\ell>0$, we denote by $\mathcal{F}_\ell(x)$ the $\ell$-neighborhood of $x$ in the leaf $\mathcal{F}(x)$ under its intrinsic topology.
Given a partially hyperbolic diffeomorphism, we simply call a plaque family corresponding to the splitting $E^c\oplus E^u$ as a \emph{$cu$-plaque family}.

 Given a partially hyperbolic diffeomorphism $f$ with   $\dim(E^c)=1$ and $\chi>0$, define the Pesin block of the following form:
 \[{\rm NUH}_\chi=\big\{x\in M: \|Df^{-n}|_{E^c(x)}\|\leq e^{-n\chi} \textrm{~for any $n\in\mathbb{N}$} \big\}. \] 
 \begin{Theorem}\label{thm.key-theorem}
 	Let $f\in\diff^1(M)$ be a partially hyperbolic diffeomorphism with  $\dim(E^c)=1$. 
 Suppose that 
 \begin{itemize}
 	\item  there exists $\chi_0>0$ such that $\int\log\|Df|_{E^c}\|\ud\mu>\chi_0$ for any $\mu\in G^u(f)$; 
 	\item 
 	let $\cW^{cu}$ be a $cu$-plaque family and  $\e_0>0$  be  given by Lemma~\ref{l.uniform-unstale-manifold} corresponding to $\chi_0$, then 
 	there exists a constant $\ell_0>0$ such that the strong stable foliation is $(\ell_0,\e_0/2)$-dense with respect to $\cW^{cu}$. 
 \end{itemize}
 Then there exist a constant $\rho_1>0$ and a hyperbolic periodic point $q$ with positive center Lyapunov exponent and whose unstable manifold has inner radius $\e_0$ 
 such that the  following properties hold. For any ergodic measure $\nu$ with $\chi_0/24<\chi^c(\nu)\leq 0$ and any $\e>0$, there exists $\xi_1>0$ such that for any  $\xi\in(0,\xi_1)$, there exist  periodic points $p_1,\cdots, p_m\in {\rm NUH}_{2(|\chi^c(\nu)|+\e)}$ of the same period $l\in\mathbb{N}$ satisfying  that 
 \begin{itemize}
 	\item each $p_i$ is homoclinically related to $q$;
 	\item $\ud(p_i,p_j)<\xi/16$ for $1\leq i, j\leq m$;
 	\item $\{p_1,\cdots,p_m\}$ is a $(l,\xi)$-separated set;
 	\item $\frac{\log m}{l}> \frac{h_\nu(f)-\e}{1+\rho_1(|\chi^c(\nu)|+\e)};$
 	\item $\ud(\delta_{\cO_{p_j}},\nu)<\rho_1(|\chi^c(\nu)|+\e)$ for any $1\leq j\leq m$. 
 \end{itemize}
%
 %
 \end{Theorem}

 The periodic orbits obtained in Theorem~\ref{thm.key-theorem} would be used to find horseshoes approaching the ergodic measure $\nu$ with $\chi^c(\nu)\leq 0$ in weak$*$-topology  and in entropy. To find periodic orbits, we first need to find some periodic pseudo-orbit which are close to the ergodic measure in weak$*$-topology and have weak (but not too weak) center Lyapunov exponent, then Liao's shadowing lemma is applied.  
 Some similar ideas can be found in \cite{BZ,YZ}, which assume  both minimality of strong foliations and  existence of $cs$-blenders and $cu$-blenders (one can refer to \cite{BD} for the definition of blenders). 
 
 \subsection{Generating orbit segments with weak hyperbolicity}

 The following is a key result to find pseudo orbits with some weak expansion along the center and satisfying the assumption in Liao's shadowing lemma.
 	\begin{Proposition}\label{p.key}
 	Let $f\in\diff^1(M)$ be a partially hyperbolic diffeomorphism with      $\dim(E^c)=1$.
 	Suppose that 
 	\begin{itemize}
 		\item  there exists $\chi_0>0$ such that $\int\log\|Df|_{E^c}\|\ud\mu>\chi_0$ for any $\mu\in G^u(f)$; 
 		\item let $\cW^{cu}$ be a $cu$-plaque family and  $\e_0>0$  be  given by Lemma~\ref{l.uniform-unstale-manifold} corresponding to $\chi_0$, then 
 		there exists a constant $\ell_0>0$ such that the strong stable foliation is $(\ell_0,\e_0/2)$-dense with respect to $\cW^{cu}$. 
 		\end{itemize}
 Then there exist  a hyperbolic periodic point $q$  and a constant $\rho_0>0$  with the following properties:
 \begin{itemize}
 	\item $\|Df^{-i}|_{E^c(q)}\|< e^{-i\chi_0}$ for any $i\in\NN$; 
	\item  for any  $C>1$, $\chi\in (-\chi_0/24,\chi_0/24)$,  $\e\in(0,\chi_0/24)$   and  $d>0$, there exist   integers $\tau=\tau(\e,d)<T=T(C,\chi,\e,d)$ such that for any $n\geq T$ and any $x\in M$ 
	satisfying \[C^{-1}\cdot e^{k(\chi-\e)}\leq  \|Df^k|_{E^c(x)}\|\leq C\cdot e^{k(\chi+\e)} \: \textrm{~for any $0\leq k\leq n,$} \]
	 there exist a point  $w_x\in W^u_{d/2}(q)$  and an integer $l_x\in\big(n,\: n+\rho_0\cdot(|\chi|+\e)\cdot n\big)$ such that 
 \begin{itemize}
 	\item  $f^{l_x}(w_x)\in\cF^s_{d/2}(q)$;
 	\item  $\ud(f^{\tau+i}(w_x),f^{i}(x))<d$ for any $i\in [\tau,n]$;
 	\item  $\|Df^{-i}|_{E^c(f^{l_x}(w_x))}\|\leq e^{-i2(|\chi|+\e)}$ for any $i\in[1,l_x]$.
 \end{itemize}
 \end{itemize}
%
 		\end{Proposition}
 	
 The idea of proving Proposition~\ref{p.key}	is that we first replace $x$ by a point $x^s$ on the local unstable manifold of a hyperbolic periodic saddle $q$, then we choose a thin $cu$-strip at $x^s$ which is $C^1$-foliated by discs tangent to a strong unstable cone field  and whose size in the center direction is exponentially small; then following the forward orbit of $x^s$, we iterate this strip  for a long time $n$ and we are interested in those points in the $cu$-strip whose forward orbits stay close to the one of $x^s$ until time $n$. Then we iterate the strip at $f^n(x^s)$ by some uniform finite time $\tau$ to make the size in the strong unstable direction at some uniform scale  so that we can apply the large deviation property for the set $G^u(f)$ given by Theorem~\ref{thm.large-deviation}.  Combining with Pliss lemma, we can find a point has "Pliss-property" which guarantees that after iterating the strip at $f^{n+\tau}(x^s)$ by $m$-times, the $cu$-strip has size at scale $\e_0$, and thus it intersects the local strong stable manifold of $q$. Our argument can guarantee that $m/n$ is at   scale-$\rho_0\cdot(|\chi|+\e)$ for some uniform constant $\rho_0.$ 
 
   In this section, a \emph{center curve} means a $C^1$-curve everywhere tangent to the center bundle. 
 
 	\proof[Proof of Proposition~\ref{p.key}] 
 	Define 
 	\[\chi_{\rm max}=\sup_{y\in M}|\log\|Df|_{E^c(y)}\||>0.\]
 	By Lemma~\ref{l.the-set-Gu(f)}, the set $G^u(f)$ is non-empty and compact. Let us define 
 	\[\chi^u=\inf\big\{ \int\log\|Df|_{E^c}\|\ud\mu\,|\, ~\mu\in G^u(f)\big\}. \] 
 	By assumption, 
 	one has 
 	$0<\chi_0<\chi^u\leq \chi_{\rm max}.$ By the local invariant property of the plaque family $\cW^{cu}$, for $\e_0$ given in the assumption,  there exists $\delta_0>0$ such that 
 	\begin{equation}\label{eq:choice-of-delta-0}
 		f^{-1}(\cW_{\delta_0}^{cu}(y))\subset\cW^{cu}_{\e_0}(f^{-1}(y))\,\textrm{~ for any $y\in M.$} 
 		\end{equation}

 		Since hyperbolic ergodic measures are approximated by periodic measures \cite{C, G2, Ge}, by the continuity of the function $\log\|Df|_{E^c}\|$, there exists a periodic point $q$ with $\chi^c(\delta_{\cO_q})>\chi_0.$ By   Pliss lemma (Lemma~\ref{l.pliss}), up to replacing $q$ by some $f^{i}(q)$, one can assume that   
 \begin{equation}\label{eq:pliss-at-q}
 		\|Df^{-i}|_{E^c(q)}\|< e^{-i\chi_0} \textrm{~ for any $i\in\mathbb{N}.$} 
 	\end{equation}

 		In the following, for simplicity, we shall assume that $q$ is a fixed point.  Up to changing a metric, we shall assume that $E^c$ is orthogonal to $E^u$.
 		By Lemma~\ref{l.uniform-unstale-manifold}, the   point $q$  has unstable manifold  $W^u_{\e_0}(q)=\cW^{cu}_{\e_0}(q)$ of size  $\e_0$ and tangent everywhere  to $E^c\oplus E^u$.  
 		\begin{Claim-numbered}\label{cl.precise-size-at-q}
 			For any $z\in W^u_{\e_0}(q)$, one has that  $\cW^{cu}_{\e_0}(z)\subset W^u(q)$.
 			\end{Claim-numbered}
 		\proof 
 		Take $z\in W^u_{\e_0}(q)=\cW^{cu}_{\e_0}(q)$. By Equation~\eqref{eq:pliss-at-q} and the second item of Lemma~\ref{l.uniform-unstale-manifold}, one has 
 			\[ \|Df^{-i}|_{E^c(z)}\|< e^{-3i\chi_0/4} \textrm{~ for any $i\in\mathbb{N}.$}\]
 		Recall that $\e_0$ is given by Lemma~\ref{l.uniform-unstale-manifold} corresponding to $\chi_0$. 
 		By the first item of Lemma~\ref{l.uniform-unstale-manifold}, for any $w\in \cW^{cu}_{\e_0}(z)$, one has   
 		\[ \|Df^{-1}|_{T_w \cW^{cu}_{\e_0}(z)}\|< e^{-\chi_0/2} .\]
 		By the local invariance of the plaque family $\cW^{cu}$, one has $f^{-1}(\cW^{cu}_{\e_0}(z))\subset  \cW^{cu}_{\e_0\cdot e^{-\chi_0/2}}(f^{-1}(z)).$
 		Once again using  the first item of Lemma~\ref{l.uniform-unstale-manifold}, one has 
 		\[ \|Df^{-2}|_{T_{w}\cW^{cu}_{\e_0}(z)}\|< e^{\chi_0/2}\cdot \|Df^{-2}|_{E^c(z)}\|< e^{-\chi_0}. \]
 		By the local invariance of the plaque familly, one has   $f^{-2}(\cW^{cu}_{\e_0}(z))\in \cW^{cu}_{\e_0\cdot e^{-\chi_0}}(f^{-2}(z)).$
 		Inductively using the arguments above, one deduces that $f^{-i}(\cW^{cu}_{\e_0}(z))\subset  \cW^{cu}_{\e_0\cdot e^{-i\chi_0/2}}(f^{-i}(z))$ for any $i\in\mathbb{N}.$ Since $z\in W^u_{\e_0}(q)$,    one   deduces that $\cW^{cu}_{\e_0}(z)\subset W^u(q)$.  
 		\endproof
 		

 	Now apply $\varphi(x)=\log\|Df|_{E^c(x)}\|$ and $\kappa=(\chi^u-\chi_0)/2$ to Theorem~\ref{thm.large-deviation}, and one obtains positive constants $r_0,a,b>0$ and a $Df$-invariant cone field $\widetilde{\mathcal{C}}^u$ with posited properties. We take a smaller $Df$-invariant cone field $\mathcal{C}^u\subset \widetilde{\mathcal{C}}^u$ such that  $f$ is uniformly expanding along $\mathcal{C}^u$, that is, there exists $\lambda^u>1$ such that
 	\[ \|Df(v)\|\geq \lambda^u, \,\textrm{~for any $x\in M$ and any unit vector $v\in \mathcal{C}^u(x)$}. \]
 	Fix a $C^1$-foliation  $\widehat{F}^u$ on the local unstable manifold $W^u_{\e_0}(q)$,  whose leaves are discs tangent to the cone field $\mathcal{C}^u$ (of course, this foliation is not $f$-invariant and not unique, and we just fix one). Let $C_u>1$ be an upper bound of the Lipschitz constant for the holonomy maps (from center curves to center curves) of this foliation $\widehat{F}^u$.

 By the uniform continuity of the center bundle, for $\e>0$, there exists $\delta\in(0,\min\{\e_0/2,\delta_0\})$ such that 
 \begin{equation}\label{eq:choice-of-delta} 
 	\big|\log\|Df|_{E^c(y_1)}\|-\log\|Df|_{E^c(y_2)}\|\big| <\e/4,   \textrm{~for any $y_1,y_2\in M$ with $\ud(y_1,y_2)<\delta$.} 
 	\end{equation}
 For any $d>0$, up to shrinking $\delta$, one can assume that $\delta<d/2$. As $f$ is  uniformly contracting and    expanding  along  $E^s$ and $\mathcal{C}^u$ respectively, and as $q$ has local unstable manifold of size $\e_0$, there exists an integer $\tau=\tau(\e,d)>0$ such that for any $k\geq\tau$ one has 
 \begin{itemize}
 	\item $f^k(\cF^s_{\ell_0}(y))\subset \cF^s_{ \delta  }(f^k(y))$ for any $y\in M$;
 	\item  
 	for disc $D$ tangent to $\mathcal{C}^u$ and of inner radius $\delta/2$,  the disc $f^k(D)$ is tangent to $\mathcal{C}^u$ and contains a disc of inner radius $r_0$;
 	\item  $f^{-k}(W^u_{\e_0}(q))\subset W^u_{\delta}(q)$.
 \end{itemize} 
%
 

 Take $n>\tau$ large (which will be precised later). Now for every     $x\in M$  satisfying 
  \[C^{-1}\cdot e^{k(\chi-\e)}\leq  \|Df^k|_{E^c(x)}\|\leq C\cdot e^{k(\chi+\e)}\: \textrm{~for any $0\leq k\leq n,$} \]  let $x^s\in\cF^{s}_{\ell_0}(x)\cap W^u_{\e_0/2}(q)$, then   one has 
  \begin{equation}\label{eq:center-iterate-of-xs}
  	C^{-1}\cdot e^{-3\tau\chi_{\rm max}}\cdot e^{i(\chi-5\e/4)}\leq \|Df^i|_{E^c(x^s)}\| \leq C\cdot e^{3\tau\chi_{\rm max}}\cdot e^{i(\chi+5\e/4)}\: \textrm{~for any $0\leq i\leq n$}.
  	\end{equation}

 Now, consider a center curve $\sigma_n\subset W^u_{\e_0}(q)$ centered at $x^s$ of length $\delta\cdot e^{-2n(|\chi|+\e)}$ and the $C^1$-foliation $\widehat{F}^u$ on $W^u_{\e_0}(q)$ whose leaves are disks tangent to the cone field $\mathcal{C}^u$. Let us  denote by 
\[\mathfrak{F}_n^{cu}(x^s):=\cup_{z\in\sigma_n}\widehat F^u_{\delta}(z) \]
which is a $C^1$-submanifold of $W^u_{\e_0}(q)$. By Claim~\ref{cl.precise-size-at-q}, one also has $\mathfrak{F}_n^{cu}(x^s)\subset \cW^{cu}_{\e_0}(x^s)\subset W^u(q).$
	We call $\cup_{z\in\partial\sigma_n}\widehat F^u_{\delta}(z) $ as the u-boundary of $\mathfrak{F}_n^{cu}(x^s)$ which has two connected components. One   defines the center-size of $\mathfrak{F}_n^{cu}(x^s)$ by the infimum of the  length of the center curves in $\mathfrak{F}_n^{cu}(x^s)$ which joins the two connected components of the u-boundary of $\mathfrak{F}_n^{cu}(x^s)$. Then the center-size of $\mathfrak{F}_n^{cu}(x^s)$ is bounded from below by $C_u^{-1}\cdot \delta\cdot e^{-2n(|\chi|+\e))}$ and from above by $C_u \cdot \delta\cdot e^{-2n(|\chi|+\e)}$.
For each $i\in[0,n]$, let $\mathfrak{F}_n^{cu}(f^i(x^s))$ be the connected component of $f^i(\mathfrak{F}_n^{cu}(x^s))\cap B(f^i(x^s),\delta)$ containing $f^i(x^s)$. By the locally invariant property   of the plaque family, one deduces that  $\mathfrak{F}_n^{cu}(f^i(x^s))\subset \cW^{cu}_\delta(f^i(x^s)).$ Besides,  one can analogously define the u-boundary and the center-size of $\mathfrak{F}_n^{cu}(f^i(x^s))$ respectively.

Let $N_1=N_1(\e,d,C)$ be an integer such that 
\begin{equation}\label{eq:choice-of-N-bounding-C}
 	C_u\cdot C\cdot e^{5\tau\chi_{\rm max}}<e^{k\e/4}\, \textrm{~and~}\, C\cdot e^{6\tau\chi_{\rm max}}<e^{k\e^2/(2\chi_{\rm max})}\, \textrm{~for any $k\geq N_1$} . 
\end{equation}
Therefore for   $n\geq N_1$, the center size of $\mathfrak{F}_n^{cu}(x^s)$ is much smaller than $\delta$. 
\begin{Claim-numbered}\label{cl.center-size-of-strip}
For any $n\geq N_1$ and  each $i\in[0,n]$, the disc $\mathfrak{F}_n^{cu}(f^i(x^s))$ has its center-size in between 	$\big(\delta\cdot e^{-4n(|\chi|+\e)},\: \delta\cdot e^{-n(|\chi|+\e/4)}\big)$. 
\end{Claim-numbered} 
\proof

By Equation~\eqref{eq:choice-of-delta},  for each $y\in \mathfrak{F}_n^{cu}(f^i(x^s))$ and any $j\leq i+1$, one has 
  \[\|Df^j|_{E^c(x^s)}\|\cdot e^{-j\e/4}\leq \|Df^j|_{E^c(y)}\|\leq \|Df^j|_{E^c(x^s)}\|\cdot e^{j\e/4},\] 
  and thus by Equation~\eqref{eq:center-iterate-of-xs}, one has  
  \begin{equation}\label{eq:center-iterate-of-cu-strip}
  	C^{-1}\cdot e^{-3\tau\chi_{\rm max}}\cdot e^{j(\chi- 3\e/2)}\leq  \|Df^j|_{E^c(y)}\|\leq C\cdot e^{3\tau\chi_{\rm max}}\cdot e^{j(\chi+3\e/2)}.
  	\end{equation}

We inductively prove this claim. Assume that the conclusion holds for $i<n$, then consider a center curve $\gamma\subset f(\mathfrak{F}_n^{cu}(f^i(x^s)))$ which joins the two  u-boundary components  of the $cu$-strip $f(\mathfrak{F}_n^{cu}(f^i(x^s)))$, then by the uniform expansion of $f $ along the cone field $\mathcal{C}^u$ and the invariance of the center bundle,  $f^{-j}(\gamma)$ is a center curve joining the two u-boundary components of $\mathfrak{F}_n^{cu}(f^{i+1-j}(x^s))$. Let $\tilde{\gamma}:[0,1]\to M$ denote the center curve $f^{-(i+1)}(\gamma)$, then 
 the length $\ell(\tilde\gamma)$ of $\tilde\gamma$ has the following bounds
\[C_u^{-1} \cdot\delta\cdot e^{-2n(|\chi|+\e)}\leq \ell(\tilde\gamma)\leq C_u\cdot \delta\cdot e^{-2n(|\chi|+\e)},\]   
and the length $\ell(\gamma)$ of the center curve $\gamma$ is given by 
\[\ell(\gamma)=\int_0^1 \|Df^{i+1}\frac{\ud\tilde\gamma(t)}{\ud t}\|\ud t 
= \int_0^1 \|Df^{i+1}|_{E^c(\tilde{\gamma}(t))}\|\cdot \|\frac{\ud\tilde\gamma(t)}{\ud t}\|\ud t. \]
By Equation~\eqref{eq:center-iterate-of-cu-strip}, one deduces that 
\[ C^{-1}\cdot e^{-3\tau\chi_{\rm max}}\cdot e^{{(i+1)}(\chi-3\e/2)}\cdot\ell(\tilde\gamma)
\leq
\ell(\gamma)
\leq C\cdot e^{3\tau\chi_{\rm max}}\cdot e^{{(i+1)}(\chi+3\e/2)}\cdot  \ell(\tilde\gamma) . \]
Combining with  Equation~\eqref{eq:choice-of-N-bounding-C}, one has 
\begin{align*}
\ell(\gamma)&\leq \delta\cdot C_u\cdot C\cdot e^{3\tau\chi_{\rm max}}\cdot e^{{(i+1)}(\chi+3\e/2)}\cdot e^{-2n(|\chi|+\e)}
\\
&<\delta\cdot e^{n\e/4}\cdot e^{n(|\chi|+3\e/2)}\cdot e^{-2n(|\chi|+\e)}
\\
&=\delta\cdot e^{-n(|\chi|+\e/4)}
\end{align*}
and 
\begin{align*}
	\ell(\gamma)&\geq \delta\cdot C_u^{-1}\cdot C^{-1}\cdot e^{-3\tau\chi_{\rm max}}\cdot e^{{(i+1)}(\chi-3\e/2)}\cdot e^{-2n(|\chi|+\e)}
	\\
	&>\delta\cdot e^{-n\e/4}\cdot e^{-n(|\chi|+3\e/2)}\cdot e^{-2n(|\chi|+\e)}
	\\
	&>\delta\cdot e^{-4n(|\chi|+\e)}.
\end{align*}
\endproof
By the choice of $\tau$, the cu-disc $f^{\tau}(\mathfrak{F}_n^{cu}(f^n(x^s)))$ is $C^1$-foliated  by discs tangent to the cone field $\mathcal{C}^u$ of radius no less than $r_0$ and its center size is bounded from below by $\delta\cdot  e^{-4n(|\chi|+\e)}\cdot e^{-\tau\chi_{\rm max}}$. 
Let $D\subset f^{\tau}(\mathfrak{F}_n^{cu}(f^n(x^s)))$ be the $C^1$-disc through the point $f^{n+\tau}(x^s)$ of radius $r_0$, then there exists a small neighborhood $\widehat D$ of $f^{n+\tau}(x^s)$ in $D$  such that for each point $y\in \widehat D$,  the disc $f^{\tau}(\mathfrak{F}_n^{cu}(f^n(x^s)))$ contains a $cu$-disc centered at $y$ of radius  $\delta_n$, where $\delta_n:=2^{-1}\cdot\delta\cdot  e^{-4n(|\chi|+\e)}\cdot e^{-\tau\chi_{\rm max}}$. 
Now we show that the $cu$-plaque of size $\delta_n$ centered at $y\in\widehat D$ is contained in $f^{\tau}(\mathfrak{F}_n^{cu}(f^n(x^s)))$, which is essentially due to Claim~\ref{cl.precise-size-at-q} and the locally invariant property of $cu$-plaques. 
\begin{Claim-numbered}\label{cl.cu-plaque-in-F-cu}
	For any $y\in \widehat D$, one has  $\cW^{cu}_{\delta_n}(y)\subset f^{\tau}(\mathfrak{F}_n^{cu}(f^n(x^s))). $
\end{Claim-numbered}
\proof
 By the choice of $\delta_0$ in  Equation~\eqref{eq:choice-of-delta-0},  one has \[f^{-1}(\cW^{cu}_{\delta_0}(f^{-n-\tau+1}(y)))\subset \cW^{cu}_{\e_0}(f^{-n-\tau}(y)).\]
Since $\mathfrak{F}_n^{cu}(x^s)\subset W_{\e_0}^u(q)$, one has that  $f^{-n-\tau}(y)\in \mathfrak{F}_n^{cu}(x^s)\subset  W_{\e_0}^u(q)$. By Claim~\ref{cl.precise-size-at-q}, one has that $\cW^{cu}_{\e_0}(f^{-n-\tau}(y))\subset W^u(q)$. As   $ \mathfrak{F}_n^{cu}(f(x^s))$ is the connected component of $f(\mathfrak{F}_n^{cu}(x^s))\cap B_{\delta}(f(x^s))$ containing $f(x^s)$, combining with Claim~\ref{cl.center-size-of-strip} and $\delta<\delta_0$,  one deduces that $\cW^{cu}_{\delta_n}(f^{-n-\tau+1}(y)))\subset \mathfrak{F}_n^{cu}(f(x^s)).$ By  applying the argument above inductively,   one can conclude.
\endproof

  Apply the disc  $D$ to Theorem~\ref{thm.large-deviation}, and one deduces 
\[\Leb_D\big(\big\{y\in D: \ud\big(\frac{1}{k}\log\|Df^k|_{E^c(y)}\|, I(\varphi)\big)\leq (\chi^u+\chi_0)/2\big\}\big)> \Leb_D(D)-a\cdot e^{-kb}\, \textrm{~ for any $k\in\NN$},\]
where $I(\varphi)=\big[\chi^u,\,\sup \{ \int\log\|Df|_{E^c}\|\ud\mu\,| ~\mu\in G^u(f) \}\big]$ is an interval. 
As $b>0$,  there exists $k_0\in\NN $ such that 
\[\Leb_D\big(\big\{y\in \widehat D: \|Df^k|_{E^c(y)}\|\geq e^{k(\chi^u+\chi_0)/2} \big\}\big)>0\,  \textrm{~for any  $k\geq k_0$.}\]
Let $N_2=N_2(\chi,\e)\in\mathbb{N}$ be an integer  such that for any $n\geq N_2$, one has 
\[\frac{24n(|\chi|+\e)}{\chi_0\cdot\rho}>k_0,\]
where $\rho=\frac{\chi^u-\chi_0}{2\chi_{\rm max}-\chi^u-\chi_0}$ is given by   Pliss Lemma (Lemma~\ref{l.pliss}) with $C=\chi_{\rm max},~\chi_1=(\chi^u+\chi_0)/2~,\chi_2=\chi_0.$
Take  
\begin{equation}~\label{eq:choice-of-tildem}
\widetilde{m}=\big[\frac{12n(|\chi|+\e)}{\chi_0\cdot\rho}\big]+1>k_0.
\end{equation}
  Then there exists a point $y_0\in \widehat D$ such that $\|Df^{\widetilde m}|_{E^c(y_0)}\|\geq e^{\widetilde{m}(\chi^u+\chi_0)/2}$.  By Pliss Lemma, there exists an integer $m\in[\widetilde m\rho, \widetilde m]$ such that 
\begin{equation}\label{eq:pliss-time}
\|Df^{-i}|_{E^c(f^m(y_0))}\|\leq e^{-i\chi_0}\, \textrm{~	for every  $1\leq i\leq m$}.
\end{equation}
\begin{Claim-numbered}~\label{c.backward-iterate}
	The $cu$-plaque $\cW^{cu}_{\e_0}(f^m(y_0))$ satisfies  that
	\begin{itemize}
		\item   $f^{-m}(\cW^{cu}_{\e_0}(f^m(y_0)))\subset \cW^{cu}_{\e_0\cdot e^{-3m\chi_0/4}}(y_0)$; 
		\item $\|Df^{-i}|_{T_w \cW^{cu}_{\e_0}(f^m(y_0))}\|<e^{-3i\chi_0/4}$  \textrm{for every $w\in \cW^{cu}_{\e_0}(f^m(y_0))$ and $1\leq i\leq m$}.
	\end{itemize} 
\end{Claim-numbered}
\proof
 By the first item in Lemma~\ref{l.uniform-unstale-manifold} and Equation~\eqref{eq:pliss-time}, for any point $w\in \cW^{cu}_{\e_0}(f^m(y_0))$,  one has \[\|Df^{-1}|_{T_w \cW^{cu}_{\e_0}(f^m(y_0))}\|<e^{\chi_0/4}\cdot \|Df^{-1}|_{E^c(f^m(y_0))}\|\leq e^{-3\chi_0/4}.\]
  By the local invariance of the plaque family, one has $f^{-1}(\cW^{cu}_{\e_0}(f^m(y_0)))\subset \cW^{cu}_{\e_0\cdot e^{-3\chi_0/4}}(f^{m-1}(y_0))$. 
  Inductively apply the above arguments, and one deduces that 
\begin{itemize}
	\item $f^{-m}(\cW^{cu}_{\e_0}(f^m(y_0)))\subset \cW^{cu}_{\e_0\cdot e^{-3m\chi_0/4}}(y_0)$; 
	\item $\|Df^{-i}|_{T_w \cW^{cu}_{\e_0}(f^m(y_0))}\|<e^{-3i\chi_0/4}$  \textrm{for every $w\in \cW^{cu}_{\e_0}(f^m(y_0))$ and $1\leq i\leq m$}.
\end{itemize}
%
\endproof
Recall that  $\delta$ is determined by $\e$ (by Equation~\eqref{eq:choice-of-delta}). There exists $N_3=N_3(\e)\in\NN$ such that 
\[  2\cdot \e_0< \delta \cdot e^{k\e}\, \textrm{~for any $k\geq N_3.$ }\] 
Thus for $n> \max\{N_1,N_3\}$, as $m\geq \widetilde{m}\rho$, by Equations \eqref{eq:choice-of-N-bounding-C} and \eqref{eq:choice-of-tildem}, 
 one has  
\begin{align*}
	\e_0\cdot e^{-3m\chi_0/4}<\e_0\cdot e^{-3\widetilde{m}\rho \chi_0/4}
	<2^{-1}\cdot\delta \cdot e^{n\e-9n(|\chi|+\e)}
	 <2^{-1}\cdot \delta\cdot e^{-4n(|\chi|+\e)}\cdot e^{-\tau\chi_{\rm max}}=\delta_n.
	\end{align*}
By Claims \ref{cl.cu-plaque-in-F-cu} and \ref{c.backward-iterate}, one gets  that
\[f^{-m}(\cW^{cu}_{\e_0}(f^m(y_0)))\subset \cW^{cu}_{\e_0\cdot e^{-3m\chi_0/4}}(y_0)\subset \cW^{cu}_{\delta_n}(y_0)\subset  f^{\tau}(\mathfrak{F}_n^{cu}(f^n(x^s))).\]
Since the strong stable foliation is $(\ell_0,\e_0/2)$-dense with respect to the plaque family $\cW^{cu}$,  there exists a point
\[z_0\in\cF^s_{\ell_0}(q)\cap \cW^{cu}_{\e_0/2}(f^m(y_0))\subset \cF^s_{\ell_0}(q)\cap  f^{m+\tau}(\mathfrak{F}_n^{cu}(f^n(x^s))).\]
 Since $\mathfrak{F}_n^{cu}(f^n(x^s)))$ is tangent to $E^c\oplus E^u$, so is $\cW^{cu}_{\e_0}(f^m(y_0))$.
As $E^c$ is dominated by  $E^u$, by the second item of    Claim~\ref{c.backward-iterate}, one has that 
\begin{equation}\label{eq:pliis-at-z0}
	\|Df^{-i}|_{E^c(z_0)}\|<  -3i\chi_0/4 \textrm{~for any $1\leq i\leq m$.}
	 \end{equation}
By the choices of $z_0$ and $\tau$, one has $f^{\tau}(z_0)\subset \cF^s_{\delta}(q)$. 
By the uniform continuity of the center bundle and Equation~\eqref{eq:pliss-at-q}, there exists an integer $t=t(\e,d)>\tau$ such that 
\[\|Df^{-i}|_{E^c(f^{\tau+t}(z_0))}\|<e^{-3i\chi_0/4} \textrm {~for any $i\leq t+\tau$}. \]
 Thus one has 
\begin{equation}\label{eq:pliss-at-widehatz0}
	\|Df^{-i}|_{E^c(f^{\tau+t}(z_0))}\|<e^{-3i\chi_0/4}  \textrm {~ for every $i\leq t+\tau+m$}.
\end{equation}
Since $t$ and $\tau$ depend only on $\e$ and $d$, there exists $N_4=N_4(\e,d)\in\mathbb{N}$ such that 
\begin{equation}\label{eq:bound-tail-time}
(3\tau+t)\cdot\chi_{\rm max}< k\cdot\e \textrm{~ for any $k\geq N_4$}.
	\end{equation}
Let $T=T(C,\e,d,\chi)=\max\big\{\tau, N_1 , N_2, N_3,N_4\big\}.$ From now on, we require that $n>T.$

Let $w_x=f^{-m-2\tau-n}(z_0)$ and $l_x=m+n+3\tau+t$. Note that $w_x\in f^{-n-\tau}(\mathfrak{F}_n^{cu}(f^n(x^s)))\subset f^{-\tau}(W^u_{\e_0}(q))\subset W^u_{\delta}(q)\subset W^u_{d/2}(q)$ and 
$f^{l_x}(w_x)=f^{\tau+t}(z_0)\in \mathcal{F}^s_{\delta}(q)\subset \mathcal{F}^s_{d/2}(q)$ since $\delta<d/2$, proving the first item for $w_x$. 

Using Equations~\eqref{eq:choice-of-tildem} and ~\eqref{eq:bound-tail-time}, and the fact that $m\leq \widetilde m$, 
one has 
\[l_x<\frac{12n(|\chi|+\e)}{\chi_0\cdot\rho}+n+\frac{n\e}{\chi_{\rm max}}\leq n+n\cdot(\frac{1}{\chi_{\rm max}}+\frac{12}{\chi_0\cdot\rho})\cdot (|\chi|+\e). \]
It suffices to take 
\[\rho_0=\frac{1}{\chi_{\rm max}}+\frac{12}{\chi_0\cdot\rho} .\]

As $z_0\in \cW^{cu}_{\e_0}(f^m(y_0))\subset f^{m+\tau}(\mathfrak{F}_n^{cu}(f^n(x^s)))$, then for every $0\leq i\leq n$ one deduces that 
$d(f^{i+\tau}(w_x),f^i(x^s))<\delta.$
Since $\ud(f^{i}(x),f^{i}(x^s))<\delta$ for every $\tau\leq i\leq n$ and $\delta/2<d$, one gets that 
$d(f^{i+\tau}(w_x),f^i(x))<2\delta<d$ for every $\tau\leq i\leq n$, proving the second item for $w_x$.

It remains to show the last item for $w_x$, that is,  
\[\|Df^{j}|_{E^c(f^{l_x-j}(w_x))}\|>e^{2j\cdot(|\chi|+\e)}\,\textrm{~for any $1\leq j\leq l_x$.} \]
By Equations~\eqref{eq:choice-of-delta} and \eqref{eq:center-iterate-of-xs}, for  each $0\leq i\leq n$   one has 
\begin{equation} \label{eq:upper-bound-of-wx}
	\|Df^{i}|_{E^c(f^{\tau}(w_x))}\|\leq e^{i\e/4}\cdot \|Df^{i}|_{E^c(x^s)}\|\leq   C\cdot e^{3\tau\chi_{\rm max}}\cdot e^{i(\chi+3\e/2)}
	\end{equation}
and 
\[\|Df^{i}|_{E^c(f^{\tau}(w_x))}\|\geq e^{-i\e/4}\cdot \|Df^{i}|_{E^c(x^s)}\|\geq C^{-1}\cdot e^{-3\tau\chi_{\rm max}}\cdot e^{i(\chi-3\e/2)}.\]
Combining with Equation~\eqref{eq:pliss-at-widehatz0} and the fact that $f^{\tau+t}(z_0)= f^{l_x}(w_x)$, one has 
%
\begin{align*}
\|Df^{l_x}|_{E^c(f(w_x))}\| 
=&\|Df^{\tau }|_{E^c(w_x)}\|\cdot\|Df^{n}|_{E^c(f^{\tau} (w_x))}\|\cdot \|Df^{\tau}|_{E^c(f^{n+\tau} (w_x))}\|\cdot \|Df^{m+t+\tau}|_{E^c(f^{n+2\tau} (w_x))}\|
\\
\geq&  e^{-\tau\chi_{\rm max}}\cdot C^{-1}\cdot e^{-3\tau\chi_{\rm max}}\cdot e^{n(\chi-3\e/2)} \cdot e^{-\tau\chi_{\rm max}}\cdot e^{3(t+\tau+m)\chi_0/4}
\\
=& C^{-1}\cdot e^{-5\tau\chi_{\rm max}}\cdot e^{n(\chi-3\e/2)}  \cdot e^{3(t+\tau+m)\chi_0/4} 
\\
\text{\tiny by Equation~\eqref{eq:choice-of-N-bounding-C}}
>&   e^{n(\chi-2\e)}\cdot e^{3(t+\tau+m)\chi_0/4} .  
\end{align*}
To show that $3(t+\tau+m)\chi_0/4+n(\chi-2\e)\geq 3(m+3\tau+t+n)(|\chi|+\e)$, it suffices to prove $3m\chi_0/4-n(|\chi|+2\e)\geq 3(m+3\tau+t+n)(|\chi|+\e)$.
 As $|\chi|+\e<\chi_0/12<\chi_{\rm max}$ and $n>T\geq N_4$, by Equation~\eqref{eq:bound-tail-time}, it suffices to show that 
 $3m\chi_0/4\geq n(|\chi|+3\e) +3(m+n)(|\chi|+\e)$. Since $|\chi|+\e<\chi_0/12$, one only needs to show that 
 $m\chi_0/2\geq n(4|\chi|+6\e)$, and this inequality holds due to the fact that 
 \[m\geq \widetilde m\rho\geq \frac{12n(|\chi|+\e)}{\chi_0} .\]
 This proves that 
\[\|Df^{l_x}|_{E^c(w_x)}\|>e^{3(m+3\tau+t+n)\cdot(|\chi|+\e)}.\]
Let $i_0\in [1,l_x]$ be the biggest integer satisfying that
	\begin{equation}\label{eq:pliss-time-at-wx} \|Df^{j}|_{E^c(f^{i_0-j}(w_x))}\|>e^{2j\cdot(|\chi|+\e)}\,\textrm{~for any $1\leq j\leq i_0$}.
		\end{equation}
Then 	we claim that $i_0=l_x.$
By Pliss lemma, one has 
\[i_0\geq l_x\cdot \frac{|\chi|+\e}{\chi_{\rm max}-2(|\chi|+\e)}>l_x\cdot \frac{|\chi|+\e}{\chi_{\rm max}}>\frac{n\e}{\chi_{\rm max}}.\]
 Assume, on the contrary, that $i_0\in(n\e/\chi_{\rm max},n+2\tau]$. 
 For each $i\in [\tau, n]$,  by Equation~\eqref{eq:upper-bound-of-wx}, one has 
 \[\|Df^i|_{E^c(w_x)}\|=\|Df^\tau|_{E^c(w_x)}\|\cdot \|Df^{i-\tau}|_{E^c(f^\tau(w_x))}\|\leq  C\cdot e^{4\tau\chi_{\rm max}}\cdot e^{(i-\tau)(\chi+3\e/2)},\]
 then combining with  Equations~\eqref{eq:choice-of-N-bounding-C} and~\eqref{eq:upper-bound-of-wx}, one has 
 \begin{align*}
 	\|Df^{i_0}|_{E^c(w_x)}\|&\leq  C\cdot e^{4\tau\chi_{\rm max}}\cdot e^{(i_0-\tau)(\chi+3\e/2)}\cdot e^{2\tau\chi_{\rm max}}
 	\\
 	&<e^{n\e^2/(2\chi_{\rm max})}\cdot e^{i_0(|\chi|+3\e/2)}
 	\\
 	&<e^{ i_0\e/2}\cdot e^{i_0(|\chi|+3\e/2)}
 	\\
 	&\leq e^{i_0(|\chi|+2\e)}
 	\end{align*}
 which contradicts with $\|Df^{i_0}|_{E^c(w_x)}\|>e^{2i_0(|\chi|+\e)}.$ Thus $i_0\in(n+2\tau, l_x]$. By Equation~\eqref{eq:pliss-at-widehatz0} and $f^{\tau+t} (z_0)=f^{l_x}(w_x)$, one has 
\[\|Df^{j}|_{E^c(f^{l_x-j}(w_x))}\|>e^{2j\cdot(|\chi|+\e)}\,\textrm{~for any $n+2\tau\leq j\leq l_x$,}\]
which implies that $i_0=l_x$, since $i_0$ is the largest integer in $[1,l_x]$ satisfying Equation~\eqref{eq:pliss-time-at-wx}. The proof of Proposition~\ref{p.key} is now completed.  
 
%
%
 	\endproof
 	According to our proof the constant $\rho_0$ is given by 
 	\[\rho_0=\frac{1}{\chi_{\rm max}}+\frac{12(2\chi_{\rm max}-\chi^u-\chi_0)}{\chi_0\cdot(\chi^u-\chi_0)},\]
 	where $\chi_{\rm max}=\sup_{y\in M}|\log\|Df|_{E^c(y)}\||$ and 	$\chi^u=\inf\big\{ \int\log\|Df|_{E^c}\|\ud\mu~|\, ~\mu\in G^u(f)\big\}. $
 	\subsection{Proof of Theorem~\ref{thm.key-theorem}}
 Given a partially hyperbolic diffeomorphism $f$ with   $\dim(E^c)=1$, and constants $C>1$, $\chi\in\mathbb{R}$ and $\e>0$, 
let us  define the following type of Pesin block:
\[{\rm \mathcal{L}}_{C,\chi,\e}=\big\{x\in M|~~ C^{-1}\cdot e^{k(\chi-\e)}\leq  \|Df^k|_{E^c(x)}\|\leq C\cdot e^{k(\chi+\e)}\,\textrm{~for any $k\in\mathbb{N}$}  \big\}. \]
Thus  ${\rm \mathcal{L}}_{C,\chi,\e}$ concerns about the set of   points whose center Lyapunov exponents are $\e$-close to $\chi$.

 	Given some separated sets, the following result allows us to find a collection of periodic points which satisfy the assumptions of  Theorem~\ref{thm.existence-of-horseshoes}. We will apply Proposition~\ref{p.key} to these separated sets to find some periodic pseudo-orbits with some weak expansion along the center bundle, then we apply Liao's shadowing lemma to get periodic points.
 	 \begin{Lemma}\label{lem:intermidiate-step}
 		Let $f\in\diff^1(M)$ be a partially hyperbolic diffeomorphism with     $\dim(E^c)=1$.
 		Suppose that 
 		\begin{itemize}
 			\item  there exists $\chi_0>0$ such that $\int\log\|Df|_{E^c}\|\ud\mu>\chi_0$ for any $\mu\in G^u(f)$; 
 			\item 
 				let $\cW^{cu}$ be a $cu$-plaque family and  $\e_0>0$  be  given by Lemma~\ref{l.uniform-unstale-manifold} corresponding to $\chi_0$, then 
 			there exists a constant $\ell_0>0$ such that the strong stable foliation is $(\ell_0,\e_0/2)$-dense with respect to $\cW^{cu}$.   			
 		\end{itemize}
 		Then there exist $\rho_0>0$ and a hyperbolic periodic point $q$ with positive center Lyapunov exponent such that for any  continuous functions $\varphi_1,\cdots,\varphi_m$   on $M$,    any $\e\in(0,\chi_0/24)$, any $\chi\in(-\chi_0/24, 0]$, any $C>1$ and any $\xi>0$,  there exists $N\in\mathbb{N}$  such that for any $n>N$, any $h\geq 0$  and any $(n,2\xi)$-separated set $S\subset  {\rm \mathcal{L}}_{C,\chi,\e}$ with cardinality $\#S\geq e^{n(h-\e/2)}$, there exists a collection $P\subset {\rm NUH}_{|\chi|+\e}$ of periodic points of the same period $l\in\big(n, n+\rho_0(|\chi|+\e)n\big)$ such that 
 		\begin{itemize}
 			\item[--] $\# P> e^{n(h-\e)};$
 			\item[--] $P$ is a $(l,\xi)$-separated set;
 			\item[--] $\ud(p,p')<\xi/16$ for any $p,p'\in P$;
 			\item[--] each $p\in P$ is homoclinically related to $q$;
 			\item[--]   for any $p\in P$, there exists $x\in S$ such that for each $1\leq j\leq m$, one has
 			\[\big|\int\varphi_j\ud\delta_{\cO_p}-\int\varphi_j\ud\frac{1}{n}\sum_{i=0}^{n-1}\delta_{f^i(x)}\big|<\big(2\rho_0(|\chi|+\e)+\e\big)\|\varphi_j\|_{C^0}, \]
 			where $\|\varphi_j\|_{C^0}=\sup_{x\in M}|\varphi_j(x)|.$
 		\end{itemize}
 	\end{Lemma}
 	\proof Let $\rho_0>0$ be the constant and $q$ be the hyperbolic periodic point given by Proposition~\ref{p.key}. 
 	Let us fix   $\e\in(0,\chi_0/24)$,  $\chi\in(-\chi_0/24, 0]$,  $C>1$, $\xi>0$ and continuous functions $\varphi_1,\cdots,\varphi_m$   on $M$.
 	
 	For $\e>0$,	by the uniform continuity of the center bundle and the uniform continuity of the functions $(\varphi_j)_{1\leq j\leq m}$, there exists $\eta_1>0$ such that for any $x,y\in M$ with $\ud(x,y)<\eta_1$, one has 
 	\begin{equation}\label{eq:uniform-continuity-of-N-functions}
 		\big|\log\|Df|_{E^c(x)}\|-\log\|Df|_{E^c(y)}\|\big|<\e/2\,
 		\textrm{~and~} \,
 		\max_{1\leq j\leq m}\frac{|\varphi_j(x)-\varphi_j(y)|}{\|\varphi_j\|_{C^0}}<\e/2.
 	\end{equation} 
 	
 	 Apply $\lambda=e^{-2(|\chi|+\e)}$ and  the dominated splitting $TM=E^s\oplus (E^c\oplus E^u)$ to Theorem~\ref{l.liao-gan-shadowing}, and one obtains the constants $L>1$ and $d_0>0$ with the posited properties.
 	
 	By Lemma~\ref{l.uniform-size-of-unstable-for-single-diffeo}, there exists $\delta_0>0$ such that if a point  $x\in M$ satisfies that  \[\|Df^{-n}|_{E^c(x)}\|\leq e^{-n(|\chi|+\e)}\, \textrm{~for any $n\in\mathbb{N}$},\] then $x$ has unstable  manifold $W^u_{\delta_0}(x)$ of size $\delta_0$ and tangent everywhere to $E^c\oplus E^u.$
 	As the bundles $E^s$ and $E^c\oplus E^u$ have uniform angle, there exists  $\eta_0>0$ such that for any $x,y\in M$ with $\ud(x,y)<\eta_0$, if    $D$ is a disc which is centered at $y$, of size $\delta_0$ and tangent  everywhere to $E^c\oplus E^u$, then  $\cF_{\delta_0}^s(x)$ has non-empty transverse intersection with $D$. 
 	
 		Apply $C>1$, $\chi$, $\e>0$ and $d=\min\{\,d_0, \, \eta_0\cdot (L+1)^{-1},\,\eta_1\cdot (L+1)^{-1},\, 32^{-1}\cdot (L+1)^{-1}\xi\,\}$ to Proposition~\ref{p.key}, and one obtains integers $\tau=\tau(\e,d)$ and $T=T(C,\chi,\e,d)$ with the posited properties.
 	
 	Let $N_1=N_1(\e,d,\xi)$ be an integer such that for any $n\geq N_1$, one has  
 	\begin{itemize}
 		\item $\rho_0(|\chi|+\e)n<e^{n\e/4}$
 		\item ${4\tau} <n\e$;
 		\item the minimal number  of $(\tau,\xi)$-Bowen balls covering the whole manifold is less than $e^{n\e/4}$.
 	\end{itemize}  
 	Let $N=\max\{T, N_1 \}$. Let $n>N$ and $S$ be a $(n,2\xi)$-separated set in ${\rm \mathcal{L}}_{C,\chi,\e}$ with cardinality $\# S> e^{n(h-\e/2)}$.  Using  the pigeonhole principle,  one gets a subset  $\widetilde{S}\subset S$ such that
 	\begin{itemize}
 		\item $\widetilde{S}$ is contained in a Bowen ball $B_{\tau}(z_0,\xi)$ for some $z_0\in M$;
 		\item   $\#\widetilde{S}>e^{n(h-3\e/4)}.$  
 	\end{itemize} 
 	By Proposition~\ref{p.key}, for any  $x\in \widetilde S$,  there exist  $w_x\in W^u_{d/2}(q)$ and an integer $l_x\in\big(n,\: n+\rho_0\cdot(|\chi|+\e)\cdot n\big)$ such that 
 	\begin{itemize}
 		\item  $f^{l_x}(w_x)\in\cF^s_{d/2}(q)$;
 		\item  $\ud(f^{\tau+i}(w_x),f^{i}(x))<d$ for any $i\in [\tau,n]$;
 		\item  $\|Df^{-i}|_{E^c(f^{l_x}(w_x))}\|\leq e^{-2i(|\chi|+\e)}$ for any $i\in[1,l_x]$.
 	\end{itemize}
 	As $\rho_0(|\chi|+\e) n<e^{n\e/4}$,  by the pigeonhole principle, there exists a subset $\widehat S\subset  \widetilde{S}$ such that 
 \begin{itemize}
 	\item $l_x$ is a constant on $\widehat S$ and we denote it by $l\in\big(n,\: n+\rho_0\cdot(\chi+\e)\cdot n\big)$;
 	\item  $\#\widehat S> e^{n(h-\e)}.$	
 	\end{itemize}
 For each $x\in\widehat S$,	since $\ud(w_x,f^{l}(w_x))\leq\ud(w_x,q) +\ud(q,f^{l}(w_x))\leq d<d_0$,	by Theorem~\ref{l.liao-gan-shadowing}, there exists a periodic point $p_x$ of period $l$ such that 
 	\begin{equation}\label{eq:px-shdaows-wx}
 		\ud(f^i(p_x),f^i(w_x))<L\cdot d \textrm{ for any $0\leq i\leq l-1$}. 
 	\end{equation}
 	This gives a collection of periodic points $P:=\{\,p_x|~x\in  \widehat{S}\;\}$ of the same period $l\in\big(n,n+\rho_0(|\chi|+\e)\big)$.  For each $x\in\widehat S$, one has $\ud(p_x,q)\leq\ud(p_x,w_x)+\ud(w_x,q)<(L+1)d$. As $(L+1)d\leq  \xi/32$, one gets that $P\subset B_{\xi/32}(q)$ which implies the third item of Lemma~\ref{lem:intermidiate-step}. It also holds that   $(L+1)d<\eta_0$,  then  by the choice of $\eta_0$, the periodic point $p_x$ is homoclinically related to $q$ which proves the fourth item of Lemma~\ref{lem:intermidiate-step}.
 	Besides,  for each $x\in\widehat S$, one has 
 	\begin{equation}\label{eq:px-shadows-x}
 		\ud(f^{\tau+i}(p_x),f^{i}(x))\leq	\ud(f^{\tau+i}(p_x),f^{i}(w_x))+	\ud(f^{\tau+i}(w_x),f^{i}(x)) <(L+1)d \textrm{~ for any $i\in [\tau,n]$. }
 	\end{equation}
 	Since $(L+1)d<\eta_1,$	for each $x\in\widehat S$,  by Equations~\eqref{eq:uniform-continuity-of-N-functions} and~\eqref{eq:px-shdaows-wx}, one has 
	\begin{equation}\label{eq:pliis-at-px}
 			 \|Df^{-i}|_{E^c(p_x)}\|\leq  e^{i\e/2}\cdot \|Df^{-i}|_{E^c(f^{l}(w_x))}\|< e^{-i(|\chi|+\e)} \textrm{ for any $i\in[1,l]$}
 			 \end{equation}
 			 which implies that $P\subset {\rm NUH}_{|\chi|+\e}$.

 	Recall that $4\tau<n\e$. Now, for each $1\leq  j\leq m$, and $x\in\widehat S,$ one has  
 		\begin{align*}
 		&\big|\int \varphi_j\ud\delta_{\cO_{p_x}}- \int \varphi_j\ud\frac{1}{n}\sum_{i=0}^{n-1}\delta_{f^i(x)}\big|
 		\\
 			& =\big|\int \varphi_j\ud\frac{1}{l}\sum_{i=0}^{l-1}\delta_{f^i(p_x)}- \int \varphi_j\ud\frac{1}{n}\sum_{i=0}^{n-1}\delta_{f^i(x)}\big|
 			\\
 			&\leq \big| \frac{1}{l}\big(\sum_{i=0}^{l-1}\varphi_j(f^i(p_x))-\sum_{i=0}^{n-1}\varphi_j(f^i(x))\big)+(\frac{1}{l}-\frac{1}{n})\sum_{i=0}^{n-1}\varphi_j(f^i(x))\big|
 			\\
 			&\leq  \frac{1}{l}\sum_{i=0}^{2\tau-1}\big|\varphi_j(f^{i}(p_x))\big|+\frac{1}{l} \sum_{i=\tau}^{n-1}\big|\varphi_j(f^{i+\tau}(p_x))-\varphi_j(f^i(x))\big|
 			+\frac{1}{l}\sum_{i=n}^{l-1}\big|\varphi_j(f^{i}(p_x))\big|
 			\\
 			&\hspace{5mm}+\frac{l-n}{n\cdot l} \sum_{i=0}^{n-1}\big|\varphi_j(f^i(x))\big| 
 			\\
 	\text{\tiny by Equations~\eqref{eq:px-shadows-x} and~\eqref{eq:uniform-continuity-of-N-functions}}	
 		&	\leq   \frac{2\tau}{l}\|\varphi_j\|_{C^0}+\frac{n-\tau}{l}\cdot\frac{\e}{2}\cdot\|\varphi_j\|_{C^0}+\frac{2(l-n)}{ l}\cdot\|\varphi_j\|_{C^0} 
 			\\
 			&\leq
 			\big(\frac{\e}{2}+\frac{\e}{2}+2\rho_0(|\chi|+\e)\big)\cdot \|\varphi_j\|_{C^0}=	\big(\e+2\rho_0(|\chi|+\e)\big)\cdot \|\varphi_j\|_{C^0},
 		\end{align*} 
 		which gives the last item of Lemma~\ref{lem:intermidiate-step}.
 
 	It remains to show that  $P$ is $(l,\xi)$-separated and $\#P> e^{n(h-\e)}$.
 	  \begin{Claim-numbered}\label{cl.separation-of-periodic-points}
 		For two different points $x,x^\prime\in \widehat{S}$, there exists $2\tau\leq k<n+\tau$ such that \[\ud(f^k(p_{x}),f^k(p_{x^\prime}))>\xi.\]
 	\end{Claim-numbered}
 	\proof 
 	Since $\widehat{S}\subset S$ is also a $(n,2\xi)$-separated and is contained in $B_\tau(z_0,\xi)$, there exists $\tau\leq i<n$ such that $\ud(f^i(x),f^i(x^\prime))>2\xi$. As $2Ld<\xi$, by Equation~\eqref{eq:px-shadows-x},  one deduces that \[\ud(f^{i+\tau}(p_{x}),f^{i+\tau}(p_{x^\prime}))\geq \ud(f^i(x),f^i(x^\prime))- \ud(f^{i+\tau}(p_{x}),f^{i}(x))-\ud(f^{i+\tau}(p_{x^\prime}),f^{i}(x^\prime))>\xi.\]
 	\endproof
 	By Claim~\ref{cl.separation-of-periodic-points}, the  set of periodic points $P=\{\,p_x|~x\in  \widehat{S}\;\}$ is $(n+\tau,\xi)$-separated and   has  cardinality $\#P=\#\widehat S> e^{n(h-\e)}$, proving the first and the second items of Lemma~\ref{lem:intermidiate-step}.	 
 	\endproof
 	
 	Now, we use Lemma~\ref{lem:intermidiate-step} to prove Theorem~\ref{thm.key-theorem}. Fix  a sequence of continuous functions $(\varphi_n)$ on the manifold $M$ which forms a dense subset of the space of continuous functions on $M$ under $C^0$-topology. Let us define  the metric on the space of Borel probability measures. The distance of two Borel probability measures $\mu_1$ and $\mu_2$  is defined as 
 	\[\ud(\mu_1,\mu_2)=\sum_{n\in\NN}\frac{|\int\varphi_n\ud\mu_1-\int\varphi_n\ud\mu_2|}{2^n(1+\|\varphi_n\|_{C^0})}, \]
 	where $\|\varphi\|_{C^0}=\sup_{x\in M}|\varphi(x)|$ for any continuous function $\varphi:M\to\RR$.
 	\proof[Proof of Theorem~\ref{thm.key-theorem}] Let $\rho_0>0$ be the constant and $q$ be the hyperbolic periodic point given by Lemma~\ref{lem:intermidiate-step}. 

 	Given  $\nu\in\cM_{erg}(f)$ with $\chi_0/24< \chi^c(\nu)\leq 0$ and   $\e\in(0,\chi_0/24)$, take  $k_0\geq 1$ such that the set 
\begin{align*}
	 {\rm Basin}_{k_0,\e}(\nu):=\big\{x\in M|~~  &e^{k(\chi^c(\nu)-\e)}\leq  \|Df^k|_{E^c(x)}\|\leq  e^{k(\chi^c(\nu)+\e)} \textrm{~and~} 
	 \\
	 &\hspace{5mm}\ud(\nu,\frac{1}{k}\sum_{i=0}^{k-1}\delta_{f^i(x)})<\e/2,\: \textrm{~for any $k\geq k_0$}\big\} 
\end{align*}
 	has $\nu$-measure no less than $1/2.$ Furthermore, there exists $C>1$ such that for any  point $x\in  {\rm Basin}_{k_0,\e}(\nu)$, one has
 	\[C^{-1}\cdot e^{k(\chi^c(\nu)-\e)}\leq  \|Df^k|_{E^c(x)}\|\leq C\cdot e^{k(\chi^c(\nu)+\e)}\,\textrm{~for any $k\in\mathbb{N}$},\]
 	in formula, ${\rm Basin}_{k_0,\e}(\nu)\subset {\rm \mathcal{L}}_{C,\chi^c(\nu),\e}.$
 	For any $\xi>0$, denote by $S(n,2\xi)$  a $(n,2\xi)$-separated set of ${\rm Basin}_{k_0,\e}(\nu)$ with maximal cardinality.  By Katok's definition of metric entropy (see Section~\ref{s.definition-of-topo-metric-entropy}), one has 
 	\[\lim_{\xi\rightarrow0}\liminf_{n\rightarrow+\infty}\frac{1}{n}\log \#S(n,2\xi)\geq h_\nu(f).\]
 	Thus there exists $\xi_1>0$ such that for any $\xi\leq\xi_1$, one has 
 	\[\liminf_{n\rightarrow+\infty}\frac{1}{n}\log \#S(n,2\xi)> h_\nu(f)-\e/4.\]
 	For any $\xi<\xi_1$, there exists an integer  $N_1>k_0$  such that
 	\begin{equation*}
 	\#S(n,2\xi) > e^{n(h_\nu(f)-\e/2)} \textrm{~ for any $n\geq N_1$.}
 	\end{equation*}
 	
 	For $\e>0$, consider an integer $m=m(\e)$ such that $2^{1-m}<\e.$ 
 	Apply the continuous functions $\varphi_1,\cdots,\varphi_m$, $\e\in(0,\chi_0/24)$, $\chi=\chi^c(\nu)$, $C>1$ and $\xi<\xi_1$ to Lemma~\ref{lem:intermidiate-step}, and one obtains an integer $N$. For any $n>\max\{N,N_1 \}$, using the $(n,2\xi)$-separated set $S(n,2\xi)$, one gets  a collection $P\subset {\rm NUH}_{|\chi|+\e}$ of periodic points of the same period $l\in\big(n, n+\rho_0(|\chi^c(\nu)|+\e)n\big)$ such that 
 	\begin{itemize}
 		\item $\# P> e^{n(h_\nu(f)-\e)};$
 		\item $P$ is a $(l,\xi)$-separated set;
 		\item $\ud(p,p')<\xi/16$ for any $p,p'\in P$;
 		\item   for any $p\in P$, there exists a point $x_p\in S(n,2\xi)$ such that for each $1\leq j\leq m$, one has
 		\[\big|\int\varphi_j\ud\delta_{\cO_p}-\int\varphi_j\ud\frac{1}{n}\sum_{i=0}^{n-1}\delta_{f^i(x_p)}\big|<\big(2\rho_0(|\chi^c(\nu)|+\e)+\e\big)\|\varphi_j\|_{C^0}.\]
 	\end{itemize}
 As $S(n,2\xi)\subset {\rm Basin}_{k_0,\e}(\nu)$ and $n\geq N_1>k_0$,   one has the following estimate
 \begin{align*}
 	\ud\big(\delta_{\cO_p},\nu\big)
 	&\leq \ud\big(\delta_{\cO_p},\frac{1}{n}\sum_{i=0}^{n-1}\delta_{f^i(x_p)}\big)+\ud\big(\frac{1}{n}\sum_{i=0}^{n-1}\delta_{f^i(x_p)},\nu\big)	\\
 	&<\sum_{j=1}^\infty\frac{\big|\int \varphi_j\ud\delta_{\cO_p}- \int \varphi_j\ud\frac{1}{n}\sum_{i=0}^{n-1}\delta_{f^i(x_p)}\big|}{2^j\|\varphi_j\|_{C^0}} +\frac{\e}{2}
 	\\	
 	&\leq 2\rho_0(|\chi^c(\nu)|+\e)+\e+\sum_{j=m+1}^\infty\frac{1}{2^j}+\frac{\e}{2}
 	\\
 	& <2\rho_0(|\chi^c(\nu)|+\e)+2\e. 	
 \end{align*}	
 Besides, one has  
 \[\frac{\log\#P}{l}> \frac{n(h_\nu(f)-\e)}{n+n\rho_0(|\chi^c(\nu)|+\e)}=\frac{h_\nu(f)-\e}{1+\rho_0(|\chi^c(\nu)|+\e)}. \]
 Now, it suffices to take $\rho_1=2\rho_0+2,$ ending the proof of Theorem~\ref{thm.key-theorem}. 
  \endproof

 \section{Approximation of ergodic measures by horseshoes: Proofs of our main results} 
 In this section, we will first show that non-hyperbolic ergodic measures are approached by horseshoes in entropy and in weak$*$-topology (i.e., proving Theorem~\ref{thm.approaching-non-hyperbolic-ergodic-measure}). Then we prove the continuity of topological entropy  and the intermediate entropy property (i.e., proving  Theorem~\ref{thm.existence}). At last, we give the proof of Theorem~\ref{thm.skeleton} and Corollary~\ref{cor:conservative-case}.
 
 \proof[Proof of Theorem \ref{thm.approaching-non-hyperbolic-ergodic-measure}]
 Let $f$ be as in the assumption. Recall that  $G^u(f)$ is compact and  varies upper semi-continuously with respect to $f$ (see Lemma~\ref{l.the-set-Gu(f)}), and thus there exist a $C^1$-open neighborhood $\cU_1$ of $f$ and a constant $\chi_0>0$ such that 
 \[\int\log\|Dg|_{E_g^c}\|\ud\mu>\chi_0\,\textrm{~for any $g\in\cU_1$ and $\mu\in G^u(g)$.}\]
 By Lemma~\ref{l.uniform-unstale-manifold}, for $\chi_0>0$, there exist a $C^1$-neighborhood $\widetilde\cU_2$ of $f$, a $cu$-plaque family $\cW^{cu}_g$ relying on $g\in\widetilde\cU_2$ continuously and $\e_0>0$ satisfying the posited properties. As the strong stable foliation of $f$ is minimal,  by Lemma~\ref{l.minimal-strong-stable-to-cu-dense}, there exist a $C^1$-neighborhood $\cU_2\subset\widetilde\cU_2$ of $f$   and $\ell_0>0$ such that  the strong stable foliation of $g\in\cU_2$ is $(\ell_0,\e_0/2)$-dense with respect to $\cW^{cu}_g$. 
 Now let $\cU=\cU_1\cap \cU_2$,  then each $g\in\cU$ satisfies the assumption of Proposition~\ref{p.key}.
Let $\chi_1= \chi_0/24$.
 
 Up to shrinking $\cU$, one can assume that 
 \[ \inf\big\{\int\log\|Dg|_{E^c_g}\|\ud\mu\;|\; g\in\cU\,\textrm{~and~} \mu\in G^u(g)\big\}>\chi_0.\]
 
Fix $g\in\cU$, then  Theorem~\ref{thm.key-theorem} gives a constant $\rho_1=\rho_1(g)>0$ for $g.$

 Let $\nu$ be an ergodic measure of $g$ with $\chi_0/24<\chi^c(\nu)\leq 0.$ Take $\e>0$ small, 
 apply $\lambda=e^{-(|\chi^c(\nu)|+\e)}$ and  the dominated splitting $TM=E_g^s\oplus (E_g^c\oplus E_g^u)$ to Theorem~\ref{thm.existence-of-horseshoes}, and one obtains $\xi_0>0$ with the posited properties.
 Apply $\nu$ and $\e>0$ to Theorem~\ref{thm.key-theorem}, and one obtains $\xi_1>0$ with posited properties. 
 For $\xi<\min\{\xi_0,\xi_1\}$, by Theorem~\ref{thm.key-theorem}, there exist finitely many periodic points $p_1,\cdots, p_m\in {\rm NUH}_{|\chi^c(\nu)|+\e}$ of the same  period $l$ such that 
 \begin{itemize}
 	\item $\ud(p_i,p_j)<\xi/16$ for any $i,j\in\{1,\cdots, m\}$;
 	\item $\{p_1,\cdots,p_m\}$ is a $(l,\xi)$-separated set;
 	\item $\frac{\log m}{l}\geq \frac{h_\nu(g)-\e}{1+\rho_1(|\chi^c(\nu)|+\e)};$
 	\item $\ud(\delta_{\cO_{p_i}},\nu)<\rho_1(|\chi^c(\nu)|+\e)$ for any $i\in\{1,\cdots, m\}$;
 \end{itemize}
 	Now, we apply Theorem~\ref{thm.existence-of-horseshoes} to the set of period points $p_1,\cdots, p_m$ and obtain a horseshoe $\La$ whose center is uniformly expanding such that 
 \begin{itemize}
 	\item $h_{top}(g, {\La})\geq \frac{\log m}{l}>\frac{h_\nu(g))-\e}{1+\rho_1(|\chi^c(\nu)|+\e)};$
 	\item  $\cM_{inv}(g,\Lambda)$ is contained in the $\e$-neighborhood of \[\big\{\sum_{i=1}^mt_i\cdot \delta_{\cO_{p_i}}|~t_i\geq 0 \textrm{~and~}\sum_{i=1}^mt_i=1\big\},\]
 	and thus in the $\rho_1(|\chi^c(\nu)|+\e)+\e$-neighborhood of $\nu.$
 \end{itemize}
Let us denote  $\chi_{\rm max}(g)= \sup_{x\in M}|\log\|Dg|_{E_g^c(x)}\||$ and $\chi^u(g)=\inf_{\mu\in G^u(g)}\int\log\|Dg|_{E_g^c}\|\ud\mu.$
By the proofs of Theorem~\ref{thm.key-theorem} and Proposition~\ref{p.key}, it suffices to take 
\[\kappa=\sup_{g\in\cU}\rho_1(g)+1=\sup_{g\in\cU}2\big(\frac{1}{\chi_{\rm max}(g)}+\frac{12(2\chi_{\rm max}(g)-\chi^u(g)-\chi_0)}{\chi_0\cdot(\chi^u(g)-\chi_0)}\big)+3.\]

Now, we need the robustly minimality of the strong stable foliation which is guaranteed by Theorem~\ref{thm.robust-minimal-ss}. Up to shrinking $\cU$, one can assume that the strong stable foliation of  $g\in\cU$ is minimal, and we will show that $\{\nu\in\mathcal{M}_{erg}(g)| \chi^c(\nu,g)\geq 0 \}$ is path connected for any $g\in\cU$. Consider $\nu_1,\nu_2\in\mathcal{M}_{erg}(g)$ with non-negative center Lyapunov exponents, then there exists a sequences of periodic orbits $\{\cO_{p_n} \}_{n\in\ZZ}$  with positive center Lyapunov exponents such that $\lim_{n\rightarrow+\infty}
\delta_{\cO_{p_n}}=\nu_1$ and $\lim_{n\rightarrow-\infty}
\delta_{\cO_{p_n}}=\nu_2.$   For each $n\in\ZZ$, by the minimality of the strong stable foliation, $\cO_{p_n}$ is homoclinically related to $\cO_{p_{n+1}}$ and thus by the Birkhoff-Smale's theorem, there exists a hyperbolic  horseshoe  $\Lambda_n$ of $g$ containing $\cO_{p_n}$ and $\cO_{p_{n+1}}$. By Theorem B in \cite{Sig},  there exists a continuous  path $\alpha_n$ on $\mathcal{M}_{erg}(g,\Lambda_n)$, which is  defined on $[1-2^{-n},1-2^{-(n+1)}]$ provided that $n\geq 0$ or defined on $[2^{n}-1,2^{n+1}-1]$ provided that $n<0$, such that the two endpoints of $\alpha_n$ are $\delta_{\cO_{p_n}}$ and $\delta_{\cO_{p_{n+1}}}$ respectively; furthermore, one can require that the path $\alpha_n$  is contained in the $2\ud(\delta_{\cO_{p_n}},\delta_{\cO_{p_{n+1}}})$-neighborhood of $\delta_{\cO_{p_n}}$ (see the comments after Theorem B in \cite{Sig}). As $\lim_{|n|\rightarrow+\infty}\ud(\delta_{\cO_{p_n}},\delta_{\cO_{p_{n+1}}})=0$, 
one gets a continuous path $\alpha:[-1,1]\to\{\nu\in\mathcal{M}_{erg}(f)| \chi^c(\nu)\geq 0 \}$ defined as 
\[ \alpha(t) =
\begin{cases}
	\nu_2 & \text{if } t=-1 ,\\
	\alpha_n(t) & \text{if } t\in [2^{n}-1,2^{n+1}-1] \text{~and~} n<0,\\
\alpha_n(t) & \text{if } t\in [1-2^{n},1-2^{n+1}] \text{~and~} n\geq 0,\\
\nu_1   & \text{if } t=1.
\end{cases} \]
 \endproof
 Now, we give the proof of Theorem~\ref{thm.skeleton} which gives some descriptions on the set of points with vanishing center Lyapunov exponents. The idea  is that the entropy of $h_{top}(\mathcal{R}_g(0))$ is bounded by lower capacity entropy which can give some separated sets,  then one can apply Lemma~\ref{lem:intermidiate-step} to the separated sets and  gets some periodic orbits satisfying Theorem~\ref{thm.existence-of-horseshoes}.
\proof[Proof of Theorem~\ref{thm.skeleton}]
We take the $C^1$-open neighborhood $\cU$ of $f$ as in  the proof of Theorem \ref{thm.approaching-non-hyperbolic-ergodic-measure}. Hence each $g\in\cU$ satisfies the hypothesis of Lemma~\ref{lem:intermidiate-step} and one obtains a constant $\rho_0>0$ and a hyperbolic periodic point $q$ with the posited properties. 

Now, we fix $g\in\cU.$ Denote by 
\[\chi_{\rm max}(g)=\sup_{x\in M}\big|\log\|Dg|_{E^c_g(x)}\|\big| \]
Let $\e>0$ and $h\in[0,h_{top}(g,\mathcal{R}_g(0))].$ Take $\e'\in(0,\chi_0/24)$   such that 
\[\frac{h-\e'}{1+\rho_0\e'}>h-\e \,\textrm{~and~}\,2\e'+\big(2\rho_0\e'+\e'\big)\chi_{\rm max}(g)<\e.\]

 By the uniform continuity of the center bundle, there exists $\eta_1>0$ such that 
\begin{equation}\label{eq:choice-of-eta1-D}
	|\log\|Dg|_{E^c_g(x)}\|-\log\|Dg|_{E^c_g(y)}\||<\e', \textrm{ for any $x,y\in M$ with $d(x,y)<\eta_1.$} \end{equation}
By the compactness of the set of Borel probability measures on $M$, there exists $\eta_2>0$ such that for any  Borel probability measures $\nu_1,\nu_2$ on $M$ with $\ud(\nu_1,\nu_2)<\eta_2$, one has 
\begin{equation}\label{eq:choice-of-eta2-D}
	\big|\int\log\|Dg|_{E^c_g}\|\ud\nu_1-\int\log\|Dg|_{E^c_g}\|\ud\nu_2\big|<\e'.
	\end{equation}

Apply $\lambda=e^{-\e'}$, $\eta_2>0$ and the dominated splitting $TM=E^s_g\oplus(E^c_g\oplus E^u_g)$ to Theorem~\ref{thm.existence-of-horseshoes}, and one gets a constant $\xi_0>0$ with the posited properties. 

  For any $k\in\mathbb{N}$, define 
\[\mathcal{R}_g^k(0)=\big\{x\in\mathcal{R}_g(0)\; |~e^{-|n|\e'}\leq \|Dg^{n}|_{E^c_g(x)}\|\leq  e^{|n|\e'}, \textrm{~for any 
$n\in\mathbb{Z}$ with $|n|>k$} \big\}. \]
Then $\mathcal{R}_g^k(0)\subset \mathcal{R}_g^{k+1}(0)$ for any $k\in\mathbb{N}$ and 
$ \mathcal{R}_g(0)=\cup_{k\in\mathbb{N}}\mathcal{R}_g^k(0).$
By the second item in Proposition~\ref{p.basic-property-entropy},  one has $h_{top}(g, \mathcal{R}_g(0))=\lim_{k\rightarrow\infty}h_{top}(g, \mathcal{R}^k_g(0))$. As $0\leq h\leq h_{top}(g,\mathcal{R}_g(0))$, there exists $k_0\in\mathbb{N}$    such that $h_{top}(g, \mathcal{R}_g^{k_0}(0))\geq h-\e'/4.$  Let $C>1$ be a constant such that for any $x\in\mathcal{R}_g^{k_0}(0)$, one has 
\[C^{-1}\cdot e^{-|n|\e'}\leq \|Dg^{n}|_{E_g^c(x)}\|\leq C\cdot  e^{|n|\e'} \textrm{~for any 
	$n\in\mathbb{Z}$ .}  \]
For $\xi>0$ and $n\in\mathbb{N}$, fix  a $(n,2\xi)$-separated set $S(n,2\xi)$ of $\mathcal{R}_g^{k_0}(0)$  with maximal cardinality. By the third item of Proposition~\ref{p.basic-property-entropy}, one has $\underline{Ch}_{top}(g, \mathcal{R}_g^{k_0}(0))\geq h_{top}(g, \mathcal{R}_g^{k_0}(0))\geq h-\e'/4.$ Then there exists $\xi_1>0$ such that for any $\xi\leq \xi_1$, there exists $k_0<\widehat N=\widehat N(\xi,\e')\in\mathbb{N}$ such that  
$\#S(n,2\xi)\geq e^{n(h-\e'/2)}$  for  any $n>\widehat N$.

Now, apply the continuous function $\varphi(x)=\log\|Dg|_{E_g^c(x)}\|$, $\chi=0$, $\e'\in(0,\chi_0/24)$, $C>1$, and $0<\xi<\min\{\xi_0,\xi_1\}$ to Lemma~\ref{lem:intermidiate-step}, and one gets an integer $N\in\mathbb{N}$. 
For $n>\max\{N,\widehat N \}$, apply Lemma~\ref{lem:intermidiate-step} to the $(n,2\xi)$-separated set $S(n,2\xi)$ and one obtains a collection $P\subset {\rm NUH}_{\e'}$ of periodic points   with the same period $l\in(n, n+\rho_0\e' n)$ such that 
\begin{itemize}
	\item $\# P> e^{n(h-\e')};$
	\item $P$ is a $(l,\xi)$-separated set;
	\item $\ud(p,p')<\xi/16$ for any $p,p'\in P$;
	\item   for any $p\in P$, there exists $x_p\in S(n,2\xi)$ such that  
	\[\big|\int\log\|Dg|_{E^c_g}\|\ud \delta_{\cO_{p}}-\int\log\|Dg|_{E^c_g}\|\ud\frac{1}{n}\sum_{i=0}^{n-1}\delta_{f^i(x_p)}\big|<\big(2\rho_0\e'+\e'\big)\chi_{\rm max}(g). \] 
\end{itemize}
For any  $p\in P$, 
 as $S(n,2\xi)\subset \mathcal{R}^{k_0}_g(0)$ and  $n>k_0$,   one has 
\begin{equation}\label{eq:upper-bound-of-CLE-D}
	\begin{split}
	\int\log\|Dg|_{E^c_g}\|\ud\delta_{\cO_{p}}&<\int\log\|Dg|_{E^c_g}\|\ud\frac{1}{n}\sum_{i=0}^{n-1}\delta_{f^i(x)}+\big(2\rho_0\e'+\e'\big)\chi_{\rm max}(g) 
	\\
	&\leq\e'+\big(2\rho_0\e'+\e'\big)\chi_{\rm max}(g).
	\end{split}
	\end{equation}
As $\xi<\xi_0$, apply $P$ to  Theorem~\ref{thm.existence-of-horseshoes},  and one  gets a horseshoe $\Lambda$ whose center is uniformly expanding such that 
\begin{itemize}
	\item $h_{top}(g,\Lambda)\geq \frac{\log\#P}{l}>\frac{h-\e'}{1+\rho_0\e'}>h-\e;$
	\item each    $\mu\in\mathcal{M}_{inv}(g,\Lambda)$ belongs to   the $
	\eta_2$-neighborhood of $\{\sum_{p\in P}t_p\delta_{\cO_{p}}|\,t_p\geq 0\, ,\sum t_p=1\}.$
\end{itemize}
By Equations~\eqref{eq:choice-of-eta2-D} and~\eqref{eq:upper-bound-of-CLE-D}, for each $\mu\in\cM_{inv}(g,\Lambda)$, one has 
\[\int \|Dg|_{E_g^c}\|\ud\mu\leq 2\e'+\big(2\rho_0\e'+\e'\big)\chi_{\rm max}(g)<\e.\]
Since $\Lambda$ has expanding center, one has $\int \|Dg|_{E_g^c}\|\ud\mu>0.$ This ends the proof of Theorem~\ref{thm.skeleton}.
\endproof 
 
Now, we turn to prove the continuity of topological entropy and the intermediate entropy property. 
\proof[Proof of Theorem~\ref{thm.existence}]
 Let $\cU$ be the $C^1$-open neighborhood of $f$ given by Theorem~\ref{thm.approaching-non-hyperbolic-ergodic-measure}. For partially hyperbolic diffeomorphisms with one-dimensional center,  by Corollary C and Theorem D in \cite{LVY} (see also Theorem 1.1 in \cite{DFPV}), one has
 \begin{itemize}
 	\item   the metric entropy varies upper semi-continuously, and in particular, there exist measures of  maximal entropy;
\item the function 
 $g\in\cU\mapsto h_{top}(g)$
 is upper semi-continuous. 
  \end{itemize}
 For each $g\in\cU$, let $\nu$ be an ergodic measure  of maximal entropy. By Theorem~\ref{thm.approaching-non-hyperbolic-ergodic-measure} and Katok's result \cite{Ka} (see also \cite{Ge}), $\nu$ is approximated by horseshoes in weak$*$-topology and in entropy.  By the structural stability of horseshoes, one deduces that the topological entropy is lower semi-continuous at $g$, and hence is continuous at $g$, proving that $g\in\cU\mapsto h_{top}(g)$ is continuous.  It is classical that  horseshoes 
 satisfy the intermediate entropy property implying that $g$  also satisfies the  intermediate entropy property.
 \endproof

 At last, we give the proof of Corollary \ref{cor:conservative-case}.
 	\proof[Proof of Corollary~\ref{cor:conservative-case}]
 	By Theorem~\ref{thm.citerion-mostly-expanding}, $f$ is mostly expanding, and thus all the ergodic measures in $G^u(f)$ have positive center Lyapunov exponents. Then Corollary~\ref{cor:conservative-case} follows directly from Theorem~\ref{thm.existence}.
 	\endproof

 	%
%
%
%
%

%
%
%
%
%
%

 \bibliographystyle{plain}

\begin{thebibliography}{99}
\bibitem[ABC]{ABC} F. Abdenur, C. Bonatti and S. Crovisier, Nonuniform hyperbolicity for
$C^1$-generic diffeomorphisms, \emph{Israel J. Math.} {\bf183} (2011), 1--60.	
 	
 	
\bibitem[ABV]{ABV}J. Alves, C. Bonatti and M. Viana,  SRB measures for partially hyperbolic systems whose central direction is mostly expanding.	\emph{Invent. Math.} {\bf140} (2000), no. 2, 351--398.

\bibitem[ADLP]{AlDLV}J. Alves, C. Dias, S. Luzzatto,and V. Pinheiro, SRB measures for partially hyperbolic systems whose central direction is weakly expanding.
\emph{J. Eur. Math. Soc.}  {\bf19} (2017), no. 10, 2911--2946.

\bibitem[ALi]{AlLi}J. Alves and X. Li, 
Gibbs-Markov-Young structures with (stretched) exponential tail for partially hyperbolic attractors. \emph{Adv. Math.} {\bf279}(2015), 405--437.
 	
\bibitem[An]{A}M. Andersson,  Robust ergodic properties in partially hyperbolic dynamics. \emph{Trans. Amer. Math. Soc.} {\bf362} (2010), no. 4, 1831--1867.

\bibitem[AnV1]{AV}M.	Andersson  and  C. H. V\'asquez,  	On mostly expanding diffeomorphisms.  	\emph{Ergodic Theory Dynam. Systems} {\bf38} (2018), no. 8, 2838--2859.	

\bibitem[AnV2]{AV2}M.	Andersson  and  C. H. V\'asquez,  Statistical stability of mostly expanding diffeomorphisms.
\emph{Ann. Inst. H. Poincar\'e C Anal. Non Lin\'eaire} {\bf37} (2020), no. 6, 1245--1270.




\bibitem[BD]{BD}C. Bonatti and L. J. D\'iaz, Persistent nonhyperbolic transitive diffeomorphisms. \emph{Ann. of Math.} (2) {\bf143} (1996), no. 2, 357--396.


 
 
\bibitem[BGP]{BGP} C. Bonatti, A. Gogolev and R. Potrie,  Anomalous partially hyperbolic diffeomorphisms II: stably ergodic examples.
\emph{Invent. Math.} {\bf206} (2016), no. 3, 801--836.

\bibitem[BGHP]{BGHP}C. Bonatti, A. Gogolev, A. Hammerlindl, and R. Potrie. Anomalous partially hyperbolic diffeomorphisms III: Abundance and incoherence. \emph{Geom. Topol.}  {\bf24} (2020) no.4, 1751--1790.

\bibitem[BZ]{BZ}
Ch. Bonatti and J.~Zhang,
\newblock Periodic measures and partially hyperbolic homoclinic classes.
\newblock {\em Trans. Amer. Math. Soc.}, 372(2):755--802, 2019.


\bibitem[Bow]{B} R. Bowen, Topological entropy for noncompact sets.
\emph{Trans. Amer. Math. Soc.} {\bf184} (1973), 125--136.

%
\bibitem[Br]{Br}A. Brown,  
Smoothness of stable holonomies inside center-stable manifolds. \emph{Ergodic Theory Dynam. Systems} {\bf42} (2022), no.12, 3593--3618.


\bibitem[BCF]{BCF}J. Buzzi, S. Crovisier and T. Fisher, 
The entropy of $C^1$-diffeomorphisms without a dominated splitting. \emph{Trans. Amer. Math. Soc.} {\bf370} (2018), no.9, 6685--6734.

\bibitem[BF]{BF}J. Buzzi and T. Fisher, Entropic stability beyond partial hyperbolicity, \emph{J. Mod. Dyn.}
{\bf7} (4) (2013), 527--552.

\bibitem[BFSV]{BFSV}J. Buzzi, T. Fisher, M. Sambarino and C. V\'asquez, Maximal entropy measures
for certain partially hyperbolic, derived from Anosov systems, \emph{Ergodic Theory Dynam. Systems} 32(1) (2012), 63--79.
%

\bibitem[C]{C}S. Crovisier,   Partial hyperbolicity far from homoclinic bifurcations. \emph{Adv. Math.} {\bf226} (2011), no. 1, 673--726.



\bibitem[CPo]{CroPol}S. Crovisier and M. Poletti, \emph{Invariance principle and non-compact center foliations}. 	arXiv:2210.14989.

 
\bibitem[CPu]{CP} S. Crovisier and E. Pujals,  Essential hyperbolicity and homoclinic bifurcations: a dichotomy phenomenon/mechanism for diffeomorphisms. \emph{Invent. Math}. {\bf201}  (2015), no.2, 385--517.
 
\bibitem[CSY]{CSY}  S. Crovisier, M. Sambarino and D. Yang,  Partial hyperbolicity and homoclinic tangencies. \emph{J. Eur. Math. Soc.} {\bf17} (2015), no. 1, 1--49.
 
\bibitem[CYZ]{CYZ}S. Crovisier, D. Yang and J. Zhang,   Empirical measures of partially hyperbolic attractors. \emph{Comm. Math. Phys.} {\bf37}5 (2020), no. 1, 725--764.

 \bibitem[DFPV]{DFPV}L. J. D\'iaz, T. Fisher, M.  Pacifico and J. Vieitez, Entropy-expansiveness for partially hyperbolic diffeomorphisms. \emph{Discrete Contin. Dyn. Syst.} {\bf32} (2012), no. 12, 4195--4207.

\bibitem[DG]{DG} L.~J. D\'{\i}az and  K.~Gelfert, 
Nonhyperbolic dynamics by mingling, blending, and flip-flopping.\emph{Topology Appl.} {\bf339} (2023), Paper No. 108571, 29 pp.
 
 \bibitem[DGR]{DGR1}L.~J. D\'{\i}az, K.~Gelfert and M. Rams, Nonhyperbolic step skew-products: ergodic approximation, 
\emph{ Ann. Inst. H. Poincar\'e C Anal. Non Lin\'eaire}  {\bf34} (2017), no. 6, 1561--1598.
 

\bibitem[DGS]{DGS} L.~J. D\'{\i}az, K.~Gelfert and B.~Santiago, 
  Weak{$*$} and entropy approximation of nonhyperbolic measures: a
geometrical approach.
 {\em Math. Proc. Cambridge Philos. Soc.}, 169(3):507--545, 2020.




\bibitem[G1]{G}S. Gan,   A generalized shadowing lemma. \emph{Discrete Contin. Dyn. Syst.} {\bf8} (2002), no. 3, 627--632.

\bibitem[G2]{G2}S. Gan,   Horseshoe and entropy for $C^1$ surface diffeomorphisms. \emph{Nonlinearity} {\bf15} (2002), no. 3, 841--848.  


\bibitem[Ge]{Ge}K. Gelfert, Horseshoes for diffeomorphisms preserving hyperbolic measures, \emph{Math. Z.} {\bf283} (3-4) (2016), 685--701.


\bibitem[Go]{Go}N. Gourmelon, Adapted metrics for dominated splittings. \emph{Ergodic Theory and Dynamical Systems}, {\bf27} (2007), no.6, 1839--1849. 

\bibitem[GSW]{GSW}L. Guan, P. Sun and W. Wu, Measures of intermediate entropies and homogeneous dynamics. \emph{Nonlinearity} {\bf9} (2017), 3349--3361.






\bibitem[HPS]{HPS}M. W. Hirsch, C. C. Pugh and M. Shub,  \emph{Invariant manifolds}. Lecture Notes in Mathematics, Vol. 583. Springer-Verlag, Berlin-New York, 1977. ii+149 pp.


\bibitem[HYY]{HYY}Y. Hua, F. Yang and J.  Yang,  
A new criterion of physical measures for partially hyperbolic diffeomorphisms.
\emph{Trans. Amer. Math. Soc.} {\bf373} (2020), no. 1, 385--417.


\bibitem[Ka]{Ka}A. Katok, Lyapunov exponents, entropy and periodic orbits for diffeomorphisms, \emph{Publ. Math. Inst. Hautes \'Etudes Sci.} {\bf51} (1980), 137--173.

\bibitem[KM]{KM}A. Katok and L. Mendoza, \emph{Dynamical systems with nonuniformly hyperbolic behavior}, Supplement to   \emph{Introduction to the modern theory of dynamical systems}, by A. Katok and B.
Hasselblatt, Encyclopedia of Mathematics and its Applications, Volume 54 (Cambridge University Press, Cambridge, 1995).
 


\bibitem[Le]{Le}F. Ledrappier,  Propri\'et\'es ergodiques des mesures de Sinai. Inst. Hautes \'Etudes Sci. Publ. Math. {\bf59} (1984), 163--188.  

\bibitem[LeS]{LeS}F. Ledrappier and J-M. Strelcyn, 
A proof of the estimation from below in Pesin's entropy formula. {\em Ergodic Theory Dynam. Systems} {\bf2} (1982), no.2, 203--219.

\bibitem[LeYo1]{LeYo1}F.  Ledrappier and L. Young,  The metric entropy of diffeomorphisms. I. Characterization of measures satisfying Pesin's entropy formula. \emph{Ann. of Math.} (2) {\bf122} (1985), no. 3, 509--539.

\bibitem[LeYo2]{LeYo2}F.  Ledrappier and L. Young,  The metric entropy of diffeomorphisms. II. Relations between entropy, exponents and dimension. \emph{Ann. of Math.}  (2) {\bf122} (1985)  no. 3, 540--574.

\bibitem[LSWW]{LSWW} M. Li, Y. Shi, S. Wang and X. Wang, Measures of intermediate entropies for star vector fields, \emph{Israel J. Math.}  {\bf240} (2020), no. 2, 791--819.

\bibitem[LVY]{LVY}G. Liao, M. Viana and J. Yang, The entropy conjecture for diffeomorphisms away from tangencies,  \emph{J. Eur. Math. Soc.}  {\bf15} (2013), no. 6, 2043--2060.

\bibitem[Li]{L}S. T. Liao, An existence theorem for periodic orbits, \emph{Acta Sci. Natur. Univ. Pekinensis}, {\bf1} (1979), 1--20.


\bibitem[Mi]{Mi}M. Misiurewicz,  Diffeomorphism without any measure with maximal entropy. \emph{Bull. Acad. Polon. Sci. S\'er. Sci. Math. Astronom. Phys.} {\bf21}(1973), 903--910.

\bibitem[N]{N}S. Newhouse, Continuity properties of entropy. \emph{Ann. of Math.} (2) {\bf129} (1989), no. 2, 215--235.

\bibitem[Pe]{Pe}Y. Pesin,  
\emph{Dimension theory in dynamical systems.} 
Contemporary views and applications.  Chicago Lectures in Math.
University of Chicago Press, Chicago, IL, 1997. xii+304 pp.

\bibitem[PeSi]{PS}Y. Pesin and Y. Sinai,   Gibbs measures for partially hyperbolic attractors. \emph{Ergodic Theory Dynam. Systems} {\bf2} (1982), no. 3-4, 417--438.


\bibitem[Pl]{Pl}V. Pliss, On a conjecture due to Smale, \emph{Differ. Uravn.} {\bf8} (1972), 262--268.



\bibitem[SY]{SY}R. Saghin and J. Yang, Continuity of topological entropy for perturbation of time-one maps of hyperbolic flows, \emph{Israel J. Math.} {\bf215} (2) (2016), 857--875.

\bibitem[Sig]{Sig}K. Sigmund, On the connectedness of ergodic systems, \emph{Manuscripta Math}. {\bf 22} (1) (1977),
27--32.

\bibitem[S1]{S3}P. Sun,  Zero-entropy invariant measures for skew product diffeomorphisms.
\emph{Ergodic Theory Dynam. Systems} {\bf30} (2010), no. 3, 923--930.

\bibitem[S2]{S}P. Sun,  Measures of intermediate entropies for skew product diffeomorphisms.
\emph{Discrete Contin. Dyn. Syst}. {\bf27} (2010), no. 3, 1219--1231.

\bibitem[S3]{S2}P.  Sun,  Density of metric entropies for linear toral automorphisms. \emph{Dyn. Syst.} {\bf27} (2012), no. 2, 197--204.

\bibitem[SW]{SW}W. Sun and  Z. Wang, Lyapunov exponents of hyperbolic measures and hyperbolic periodic orbits.
\emph{Trans. Amer. Math. Soc}. {\bf362} (2010), no. 8, 4267--4282.

\bibitem[TY]{TY}A. Tahzibi and J. Yang, Invariance principle and rigidity of high entropy measures.
\emph{Trans. Amer. Math. Soc.} {\bf371} (2019), no.2, 1231--1251.

\bibitem[U]{U}R. Ures,  Intrinsic ergodicity of partially hyperbolic diffeomorphisms with a hyperbolic linear part. \emph{Proc. Amer. Math. Soc.} {\bf 140} (2012), no.6, 1973--1985.

\bibitem[Y]{Y}J. Yang,  Entropy along expanding foliations. \emph{Adv. Math.} {\bf389} (2021), Paper No. 107893, 39 pp.

\bibitem[YZ]{YZ}D.~Yang and J.~Zhang,  Non-hyperbolic ergodic measures and horseshoes in partially hyperbolic homoclinic classes.  {\em J. Inst. Math. Jussieu}, {\bf19} (2020) no. 5, 1765--1792.

\bibitem[Yo]{Yo}Y. Yomdin, Volume growth and entropy, \emph{Israel J. Math.} {\bf57} no. 3 (1987), 285--300.

 \end{thebibliography}

\vspace{2mm}
 
 \begin{tabular}{l l l}
	\emph{\normalsize Jinhua Zhang}
	\medskip\\
	\small School of Mathematical Sciences, Beihang University, Beijing, 100191, P. R.  China\\ 
	\small \texttt{jinhua$\_$zhang@buaa.edu.cn}; \texttt{zjh200889@gmail.com}
	
\end{tabular}

\end{document}